\documentclass[10pt,a4paper]{article}
\usepackage{amsfonts}
\usepackage{latexsym}
\usepackage{cite}
\usepackage{amsmath,amsfonts,latexsym,amssymb}
\usepackage[mathscr]{eucal}
\usepackage{cases,color}
\usepackage{amsthm}

\usepackage[bf,small]{caption2}
\usepackage{float}
\usepackage{graphicx}
\usepackage{amsmath}
\usepackage{amssymb}
\usepackage[all]{xy}

\newtheorem{theorem}{theorem}[section]
\newtheorem{thm}[theorem]{Theorem}
\newtheorem{lem}[theorem]{Lemma}

\newtheorem{prop}[theorem]{Proposition}
\newtheorem{cor}[theorem]{Corollary}

\newtheorem{defn}[theorem]{Definition}
\newtheorem{exmp}[theorem]{Example}

\newtheorem{rmk}[theorem]{Remark}
\newtheorem{nota}[theorem]{Notation}

\begin{document}

\title{\textbf{${\rm SL}(3,\mathbb{C})$-character varieties of two-generator groups}}
\author{\Large Haimiao Chen \footnote{Email: \emph{chenhm@math.pku.edu.cn}}  \\
\normalsize \em{Beijing Technology and Business University, Beijing, China} }
\date{}
\maketitle

\begin{abstract}
  We present an efficient method for computing the irreducible ${\rm SL}(3,\mathbb{C})$-character varieties of two-generator groups.

\end{abstract}

\section{Introduction}

Given a finitely generated group $\Gamma$ and a linear group $G\le{\rm GL}(n,\mathbb{C})$, a {\it $G$-representation} of $\Gamma$ is simply a homomorphism $\rho:\Gamma\to G$, and its {\it character} is the function $\chi_\rho:\Gamma\to\mathbb{C}$ sending $x$ to ${\rm tr}(\rho(x))$. The {\it $G$-representation variety} of $\Gamma$ is $\mathcal{R}_G(\Gamma):=\hom(\Gamma,G)$, on which $G$ acts via conjugation, and the {\it $G$-character variety} $\mathcal{X}_G(\Gamma)$ is defined as the GIT quotient $\mathcal{R}_G(\Gamma)//G$.
Let $\mathcal{R}^{\rm irr}_G(\Gamma)$ denote the open subset of $\mathcal{R}_G(\Gamma)$ consisting of irreducible representations, and let
$$\mathcal{X}^{\rm irr}_G(\Gamma)=\{\chi_\rho\colon\rho\in\mathcal{R}^{\rm irr}_G(\Gamma)\}.$$
It is known that $\mathcal{X}^{\rm irr}_G(\Gamma)$ can be identified with an open subset of $\mathcal{X}_G(\Gamma)$, and the map $\rho\mapsto\chi_\rho$ factors through to give an isomorphism
$\mathcal{R}^{\rm irr}_G(\Gamma)/G\cong\mathcal{X}^{\rm irr}_G(\Gamma)$.
Passing to character variety removes some redundance; however, in some situations, it is necessary to work with representations. For more knowledge, see \cite{Si12} and the references therein.

Since the seminal works \cite{CS83,Th82}, character variety has been playing a significant role in low-dimensional topology, in that, much geometric and topological information on 3-manifolds $M$ is carried by $\mathcal{X}_G(\pi_1(M))$. For instance, CR and complex hyperbolic structure involve ${\rm SU}(2,1)\le{\rm SL}(3,\mathbb{C})$ \cite{De15,Pa12,PW17}.

When $n\ge 3$, as pointed out by \cite{FW17}, a main concern is to collect a sufficient amount of examples which will give directions for further research. However, in general it seems to be difficult to do explicit computations. Even for $n=3$, till now few results on representation/character varieties have been seen in the literature, except for \cite{DF15,FGKR16,FKR15,GW18,He14,HMP16,MP16,MS16}.

In this paper, we present a method to determine $\mathcal{X}^{\rm irr}_{{\rm SL}(3,\mathbb{C})}(\Gamma)$ for any two-generator group $\Gamma$, (the reducible part of $\mathcal{X}_{{\rm SL}(3,\mathbb{C})}(\Gamma)$ is, as usual, relatively easy to study). A representation $\Gamma\to{\rm SL}(3,\mathbb{C})$ is the same as a pair $(\mathbf{a},\mathbf{b})$ satisfying some equations $\mathbf{s}_j=\mathbf{0}$, each of which is equivalent to ${\rm tr}(\mathbf{u}\mathbf{s}_j)=0$ for all $\mathbf{u}$ in the vector space $\mathcal{M}$ of $3\times 3$ matrices. We provide bases (in Theorem \ref{thm:basis-1} and \ref{thm:basis-2}) in terms of $\mathbf{a},\mathbf{b}$ for $\mathcal{M}$, so as to reduce ${\rm tr}(\mathbf{u}\mathbf{s}_j)=0$ to equations of traces of words in $\mathbf{a},\mathbf{b}$.

A key ingredient (see Theorem \ref{thm:symmetric}) is the observation that, ${\rm tr}(\mathbf{a}\mathbf{b}\mathbf{a}^{-1}\mathbf{b}^{-1})={\rm tr}(\mathbf{b}\mathbf{a}\mathbf{b}^{-1}\mathbf{a}^{-1})$ if and only if $\mathbf{a},\mathbf{b}$ can be simultaneously conjugated to be symmetric. A remarkable result of Lawton \cite{La07} describes $\mathcal{X}_{{\rm SL}(3,\mathbb{C})}(F_2)$ as a hypersurface in $\mathbb{C}^9$ which is a double branched cover of $\mathbb{C}^8$, and the branching locus exactly consists of the characters with ${\rm tr}(\mathbf{a}\mathbf{b}\mathbf{a}^{-1}\mathbf{b}^{-1})={\rm tr}(\mathbf{b}\mathbf{a}\mathbf{b}^{-1}\mathbf{a}^{-1})$. Our result gives a geometric interpretation for these characters. When ${\rm tr}(\mathbf{a}\mathbf{b}\mathbf{a}^{-1}\mathbf{b}^{-1})\ne{\rm tr}(\mathbf{b}\mathbf{a}\mathbf{b}^{-1}\mathbf{a}^{-1})$, the situation turns out to be largely different, and the computation is much simpler.

As an application, the ${\rm SL}(3,\mathbb{C})$-character variety for each double twist link group is determined.

We expect that, results of this paper will promote the study of ${\rm SL}(3,\mathbb{C})$-character varieties for general groups
in future.

\begin{nota}
\rm For elements $x,y$ of some group, let $\overline{x}=x^{-1}$, $x\lrcorner y=xy\overline{x}$, $[x,y]=xy\overline{x}\overline{y}$; be careful that $\overline{x}\overline{y}=x^{-1}y^{-1}$, not $(xy)^{-1}$.
(In this paper, notations for conjugate complex numbers are never used.) For a word $r$ in $x,y$, let $\overleftarrow{r}$ or $r^\leftarrow$ denote the one obtained by writing $r$ backwards. Let ${\rm Cen}(x)$ denote the centralizer of $x$.

Let $\omega=e^{2\pi i/3}$, and let $\langle\omega\rangle=\{1,\omega,\omega^2\}$.

Let $\mathcal{M}$ denote the set of $3\times 3$ matrices with entries in $\mathbb{C}$, which is a 9-dimensional vector space over $\mathbb{C}$; let $\mathcal{M}^0$ and $\mathcal{M}^{1}$ respectively denote the subspace consisting of symmetric and skew-symmetric matrices. Let $\mathbf{e}$ denote the identity matrix.
For $\mathbf{z}\in\mathcal{M}$, let $\mathbf{z}^{\rm tr}$ denote its transpose; let
$\mathbf{z}^\neg=[z_{12},z_{13},z_{23}]^{\rm tr}$.

Given elements $\mathfrak{a}_1,\ldots,\mathfrak{a}_k$ of some vector space over $\mathbb{C}$, let
$\mathbb{C}\langle\mathfrak{a}_1,\ldots,\mathfrak{a}_k\rangle$ denote the subspace they span.
\end{nota}


\section{On $3\times 3$ unimodular matrices}

For each $\mathbf{x}\in{\rm SL}(3,\mathbb{C})$, by Cayley-Hamilton Theorem,
\begin{align}
\mathbf{x}^3-{\rm tr}(\mathbf{x})\cdot\mathbf{x}^2+{\rm tr}(\overline{\mathbf{x}})\cdot\mathbf{x}-\mathbf{e}=0. \label{eq:Hamilton-Cayley}
\end{align}
Consequently, for each $k\in\mathbb{Z}$,
\begin{align}
\mathbf{x}^k=\alpha_k(\mathbf{x})\cdot\mathbf{x}+\beta_k(\mathbf{x})\cdot\mathbf{e}
+\gamma_k(\mathbf{x})\cdot\overline{\mathbf{x}}, \label{eq:power}
\end{align}
with $\alpha_k(\mathbf{x}),\beta_k(\mathbf{x}),\gamma_k(\mathbf{x})$ polynomials in ${\rm tr}(\mathbf{x}),{\rm tr}(\overline{\mathbf{x}})$.
In particular,
\begin{align}
\mathbf{x}^3=({\rm tr}(\mathbf{x})^2-{\rm tr}(\overline{\mathbf{x}}))\cdot\mathbf{x}+(1-{\rm tr}(\mathbf{x}){\rm tr}(\overline{\mathbf{x}}))\cdot\mathbf{e}+{\rm tr}(\mathbf{x})\cdot\overline{\mathbf{x}}. \label{eq:cube}
\end{align}

Introduce the {\it standard forms}
\begin{align*}
&\mathbf{d}^{\lambda_1}_{\lambda_3}=\left[\begin{array}{ccc} \lambda_1 & 0 & 0 \\ 0 & \lambda_2 & 0 \\ 0 & 0 & \lambda_3 \end{array}\right], \qquad \text{where} \quad \lambda_2=(\lambda_1\lambda_3)^{-1}; \\
&\mathbf{h}_{\lambda}=\left[\begin{array}{ccc} 2\lambda & 0 & i\lambda \\ 0 & \lambda^{-2} & 0 \\ i\lambda & 0 & 0 \end{array}\right];  \\
&\mu\mathbf{k}, \qquad \text{where} \quad \mu\in\langle\omega\rangle, \quad
\mathbf{k}=\left[\begin{array}{ccc} 1 & 1 & 0 \\ 1 & 1 & i \\ 0 & i & 1 \end{array}\right].
\end{align*}
Clearly, each element of ${\rm SL}(3,\mathbb{C})$ can be conjugated into standard form.

Put
$$\mathbf{g}=\left[\begin{array}{ccc} 0 & 0 & 0 \\ 0 & 1 & 0 \\ 0 & 0 & 0 \end{array}\right], \qquad
\mathbf{f}=\left[\begin{array}{ccc} 0 & 1 & 0 \\ 1 & 0 & i \\ 0 & i & 0 \end{array}\right], \qquad
\mathbf{f}^2=\left[\begin{array}{ccc} 1 & 0 & i \\ 0 & 0 & 0 \\ i & 0 & -1 \end{array}\right].$$
We have $\mathbf{g}\mathbf{f}=\mathbf{f}\mathbf{g}=\mathbf{f}^3=\mathbf{0}$, and
\begin{align*}
\mathbf{k}=\mathbf{e}+\mathbf{f}, \qquad   \mathbf{h}_\lambda=\lambda\mathbf{e}+(\lambda^{-2}-\lambda)\mathbf{g}+\lambda\mathbf{f}^2.
\end{align*}

\begin{defn}
\rm An element $\mathbf{x}\in{\rm SL}(3,\mathbb{C})$ is called {\it ordinary} if it satisfies any one of the following conditions which are equivalent to each other:
\begin{itemize}
  \item[1)] $\mathbf{x}$ is neither conjugate to $\mathbf{d}^{\lambda}_{\lambda}$ for any $\lambda$, nor to $\mathbf{h}_{\kappa}$ with $\kappa\in\langle\omega\rangle$;
  \item[2)] the degree of the minimal polynomial of $\mathbf{x}$ is 3;
  \item[3)] each eigenspace of $\mathbf{x}$ is one-dimensional.
\end{itemize}
Call $\mathbf{x}$ {\it special} if it is not ordinary.
\end{defn}

\begin{lem} \label{lem:centralizer-ordinary}
If $\mathbf{x}\in{\rm SL}(3,\mathbb{C})$ is ordinary, then each element of ${\rm Cen}(\mathbf{x})$ is a linear combination of $\mathbf{e},\mathbf{x},\overline{\mathbf{x}}$, and has a square root.
\end{lem}

\begin{proof}
It suffices to consider standard forms.
\begin{itemize}
  \item If $\mathbf{x}=\mathbf{d}^{\lambda_1}_{\lambda_3}$ with $\lambda_1\ne\lambda_2\ne\lambda_3\ne\lambda_1$, then
        $${\rm Cen}(\mathbf{x})=\{\mathbf{d}^{\mu_1}_{\mu_3}\colon\mu_1\mu_3\ne0\}
        \subset\mathbb{C}\langle\mathbf{e},\mathbf{x},\overline{\mathbf{x}}\rangle;$$
        clearly each element has a square root.
  \item If $\mathbf{x}=\mathbf{h}_{\lambda}$ with $\lambda^3\ne 1$, then $\overline{\mathbf{x}}=\lambda^{-1}\mathbf{e}+(\lambda^2-\lambda^{-1})\mathbf{g}-\lambda^{-1}\mathbf{f}^2$, and
        $${\rm Cen}(\mathbf{x})=\big\{\mathbf{h}_{t,c}:=t\mathbf{e}+(t^{-2}-t)\mathbf{g}+c\mathbf{f}^2\colon t\ne 0\big\}\subset\mathbb{C}\langle\mathbf{e},\mathbf{x},\overline{\mathbf{x}}\rangle;$$
        each element has a square root, as $(\mathbf{h}_{t,c})^2=\mathbf{h}_{t^2,2tc}$.
  \item If $\mathbf{x}=\mu\mathbf{k}$ with $\mu^3=1$, then
        $${\rm Cen}(\mathbf{x})=\{a\mathbf{e}+b\mathbf{f}+c\mathbf{f}^2\colon a\in\langle\omega\rangle\}\subset\mathbb{C}\langle\mathbf{e},\mathbf{x},\overline{\mathbf{x}}\rangle;$$
        each element has a square root: $(a\mathbf{e}+b\mathbf{f}+c\mathbf{f}^2)^2=a^2\mathbf{e}+2ab\mathbf{f}+(b^2+2ac)\mathbf{f}^2$.
\end{itemize}
\end{proof}

Let $F_2=\langle z_1,z_2\mid - \rangle$ denote the rank two free group. Given $\mathbf{a}, \mathbf{b}\in{\rm SL}(3,\mathbb{C})$, let $\sigma_{\mathbf{a}, \mathbf{b}}: F_2\to{\rm SL}(3,\mathbb{C})$ send $z_1$, $z_2$ to $\mathbf{a}$, $\mathbf{b}$, respectively, and let $\chi_{\mathbf{a}, \mathbf{b}}$ denote its character. Each representation is of this form. We usually do not distinguish $(\mathbf{a},\mathbf{b})$ and $\sigma_{\mathbf{a}, \mathbf{b}}$.
Adopt the notations $t_{\overline{1}^22^3}={\rm tr}(\overline{\mathbf{a}}^2\mathbf{b}^3)$,
$t_{2\overline{1}^3\overline{2}1^4}={\rm tr}(\mathbf{b}\overline{\mathbf{a}}^3\overline{\mathbf{b}}\mathbf{a}^4)$,
and so forth; in particular,
\begin{align*}
t_1={\rm tr}(\mathbf{a}), \qquad  t_2={\rm tr}(\mathbf{b}), \qquad  t_{\overline{1}}={\rm tr}(\overline{\mathbf{a}}),
\qquad t_{\overline{2}}={\rm tr}(\overline{\mathbf{b}}), \\
t_{12}={\rm tr}(\mathbf{a}\mathbf{b}), \ \ t_{\overline{1}2}={\rm tr}(\overline{\mathbf{a}}\mathbf{b}), \ \
t_{1\overline{2}}={\rm tr}(\mathbf{a}\overline{\mathbf{b}}), \quad t_{\overline{1}\overline{2}}={\rm tr}(\overline{\mathbf{a}}\overline{\mathbf{b}}), \\
t_{12\overline{1}\overline{2}}={\rm tr}([\mathbf{a},\mathbf{b}]), \qquad  t_{21\overline{2}\overline{1}}={\rm tr}([\mathbf{b},\mathbf{a}]).
\end{align*}

According to \cite{La07} Theorem 8, the ${\rm SL}(3,\mathbb{C})$-character variety of $F_2$ can be embedded in $\mathbb{C}^9$ as a hypersurface:
\begin{align}
\mathcal{X}_{{\rm SL}(3,\mathbb{C})}(F_2)&\cong\{(s_1,s_{\overline{1}},s_2,s_{\overline{2}},s_3,s_{\overline{3}}, s_4, s_{\overline{4}},s_5)\colon s_5^2-Ps_5+Q=0\}, \nonumber  \\
\chi_{\mathbf{a}, \mathbf{b}}&\mapsto (t_1,t_{\overline{1}},t_2,t_{\overline{2}},t_{12},t_{\overline{1}\overline{2}},t_{1\overline{2}}, t_{\overline{1}2}, t_{12\overline{1}\overline{2}}),  \label{eq:char-variety-F2}
\end{align}
for certain $P,Q\in \mathbb{C}[s_1,s_{\overline{1}},s_2,s_{\overline{2}},s_3,s_{\overline{3}}, s_4, s_{\overline{4}}]$, with $\deg P=4,\deg Q=6$.
Geometrically, (\ref{eq:char-variety-F2}) presents $\mathcal{X}_{{\rm SL}(3,\mathbb{C})}(F_2)$ as a branched double cover of $\mathbb{C}^8$, with branching locus defined by $P^2-4Q=0$.

The precise meaning is: for each word $\mathbf{w}$ in $\mathbf{a},\mathbf{b}$, ${\rm tr}(\mathbf{w})$ can be written as a polynomial in $t_1, t_{\overline{1}}, t_2, t_{\overline{2}}, t_{12}, t_{\overline{1}\overline{2}}, t_{1\overline{2}}, t_{\overline{1}2}$ (see Proposition \ref{prop:express}). Actually, 
\begin{align*}
t_{12\overline{12}}+t_{21\overline{2}\overline{1}}&=P(t_1, t_{\overline{1}}, t_2, t_{\overline{2}}, t_{12}, t_{\overline{1}\overline{2}}, t_{1\overline{2}}, t_{\overline{1}2}), \\ 
t_{12\overline{1}\overline{2}}\cdot t_{21\overline{2}\overline{1}}&=Q(t_1, t_{\overline{1}}, t_2, t_{\overline{2}}, t_{12}, t_{\overline{1}\overline{2}}, t_{1\overline{2}}, t_{\overline{1}2}).
\end{align*}

\begin{nota}
\rm For a word $\mathbf{w}$ in $\mathbf{a},\mathbf{b}$, let $\theta(\mathbf{w})=\mathbf{w}+\overleftarrow{\mathbf{w}}$, and $\delta(\mathbf{w})=\mathbf{w}-\overleftarrow{\mathbf{w}}$.
\end{nota}

By \cite{Pa12} Corollary 4.3 (to Proposition 4.2),
\begin{align}
\theta(\mathbf{a}\lrcorner\mathbf{b})=\ &t_1\theta(\overline{\mathbf{a}}\mathbf{b})
+t_{\overline{1}}\theta(\mathbf{a}\mathbf{b})-(1+t_1t_{\overline{1}})\mathbf{b}
+(t_{\overline{1}2}-t_{\overline{1}}t_2)\mathbf{a}+(t_{12}-t_1t_2)\overline{\mathbf{a}} \nonumber \\
&+(t_1t_{\overline{1}}t_2+t_2-t_1t_{\overline{1}2}-t_{\overline{1}}t_{12})\mathbf{e},  \label{eq:identity1}  \\
\mathbf{a}\mathbf{b}\mathbf{a}=\ &-\theta(\mathbf{a}^2\mathbf{b})+t_1\theta(\mathbf{a}\mathbf{b})+t_2\mathbf{a}^2
+(t_{12}-t_1t_2)\mathbf{a}-t_{\overline{1}}\mathbf{b}+t_{\overline{1}2}\mathbf{e}.  \label{eq:identity2}
\end{align}
Another identity which is useful in practical computation is
\begin{align}
\mathbf{a}\mathbf{b}\mathbf{a}=t_{12}\mathbf{a}-t_{\overline{12}}\overline{\mathbf{b}}
+\overline{\mathbf{b}}\overline{\mathbf{a}}\overline{\mathbf{b}}.  \label{eq:identity3}
\end{align}
It is obtained by applying (\ref{eq:Hamilton-Cayley}) to $\mathbf{a}\mathbf{b}$ and then multiplying on the right by
$\overline{\mathbf{b}}\overline{\mathbf{a}}\overline{\mathbf{b}}$.

\begin{prop}\label{prop:express}
Each word in $\mathbf{a},\mathbf{b}$ can be written as a linear combination of
$\mathbf{e}$, $\mathbf{a}$,  $\overline{\mathbf{a}}$, $\mathbf{b}$, $\overline{\mathbf{b}}$, $\mathbf{a}\mathbf{b}$, $\mathbf{b}\mathbf{a}$, $\mathbf{a}\overline{\mathbf{b}}$, $\overline{\mathbf{b}}\mathbf{a}$, $\overline{\mathbf{a}}\mathbf{b}$, $\mathbf{b}\overline{\mathbf{a}}$, $\overline{\mathbf{a}}\overline{\mathbf{b}}$, $\overline{\mathbf{b}}\overline{\mathbf{a}}$, $\mathbf{a}\lrcorner\mathbf{b}$, $\mathbf{a}\lrcorner\overline{\mathbf{b}}$, $\mathbf{b}\lrcorner\mathbf{a}$,
$\mathbf{b}\lrcorner\overline{\mathbf{a}}$, $[\mathbf{a},\mathbf{b}]$, $[\mathbf{b},\mathbf{a}]$,
with coefficients in $\mathbb{C}[t_1, t_{\overline{1}}, t_2, t_{\overline{2}}, t_{12}, t_{\overline{1}\overline{2}}, t_{1\overline{2}}, t_{\overline{1}2}, t_{12\overline{12}}]$.
\end{prop}

\begin{proof}
By (\ref{eq:power}), it is sufficient to consider the words of the form
$\mathbf{w}=\mathbf{x}_1^{\epsilon_1}\cdots\mathbf{x}_n^{\epsilon_n}$ with each $\epsilon_i\in\{\pm1\}$ and $\{\mathbf{x}_{j},\mathbf{x}_{j+1}\}=\{\mathbf{a},\mathbf{b}\}$ for $1\le j<n$.

Whenever any of $\mathbf{a},\overline{\mathbf{a}},\mathbf{b},\overline{\mathbf{b}}$ appears twice, we can apply (\ref{eq:identity2}) to write $\mathbf{w}$ as a linear combination of shorter words. Hence we may just assume $n\le4$. Then the conclusion is deduced by applying (\ref{eq:identity1}) and the identity obtained by multiplying (\ref{eq:identity1}) on the right by $\overline{\mathbf{b}}$.
\end{proof}

\begin{cor}\label{cor:one-is-special}
\rm If $\mathbf{a}$ is special, then each word in $\mathbf{a},\mathbf{b}$ can be written as a linear combination of
$\mathbf{e},\mathbf{a},\mathbf{b},\overline{\mathbf{b}},\mathbf{a}\mathbf{b},\mathbf{b}\mathbf{a},\mathbf{a}\overline{\mathbf{b}},\overline{\mathbf{b}}\mathbf{a},\mathbf{b}\lrcorner\mathbf{a}$, with coefficients in $\mathbb{C}[t_1, t_2, t_{\overline{2}}, t_{12}, t_{1\overline{2}}]$.
\end{cor}

\begin{proof}
By Proposition \ref{prop:express}, it suffices to show this for
$$\overline{\mathbf{a}}\mathbf{b},\ \ \mathbf{b}\overline{\mathbf{a}}, \ \ \overline{\mathbf{a}}\overline{\mathbf{b}}, \ \ \overline{\mathbf{b}}\overline{\mathbf{a}}, \ \ \mathbf{a}\lrcorner\mathbf{b}, \ \ \mathbf{a}\lrcorner\overline{\mathbf{b}}, \ \
[\mathbf{a},\mathbf{b}], \ \ [\mathbf{b},\mathbf{a}].$$
Since $\overline{\mathbf{a}}\in\mathbb{C}\langle\mathbf{e},\mathbf{a}\rangle$, using (\ref{eq:identity2}) we can easily show that indeed each of these can be written as a linear combination of $\mathbf{e},\mathbf{a},\mathbf{b},\overline{\mathbf{b}},\mathbf{a}\mathbf{b},\mathbf{b}\mathbf{a},\mathbf{a}\overline{\mathbf{b}},
\overline{\mathbf{b}}\mathbf{a},\mathbf{b}\lrcorner\mathbf{a}$, with coefficients in $\mathbb{C}[t_1, t_2, t_{\overline{2}}, t_{12}, t_{1\overline{2}}]$.
\end{proof}

\section{Irreducible representations of two-generator groups}

Observe that for any $\mathbf{u},\mathbf{v},\mathbf{w}\in\mathcal{M}^1$,
\begin{align}
\det[\mathbf{u}^\neg,\mathbf{v}^\neg,\mathbf{w}^\neg]={\rm tr}(\mathbf{u}\mathbf{v}\mathbf{w}).  \label{eq:det=tr}
\end{align}

\begin{lem} \label{lem:tr-tr}
If $\mathbf{a}\in\mathcal{M}^0$, then
\begin{align*}
\det[(\mathbf{b}-\mathbf{b}^{\rm tr})^\neg,(\mathbf{a}\mathbf{b}-\mathbf{b}^{\rm tr}\mathbf{a})^\neg,(\overline{\mathbf{a}}\mathbf{b}-\mathbf{b}^{\rm tr}\overline{\mathbf{a}})^\neg]
=t_{12\overline{1}\overline{2}}-t_{21\overline{2}\overline{1}}.
\end{align*}
\end{lem}
\begin{proof}
Expanding $\mathbf{q}:=(\mathbf{b}-\mathbf{b}^{\rm tr})(\mathbf{a}\mathbf{b}-\mathbf{b}^{\rm tr}\mathbf{a})(\overline{\mathbf{a}}\mathbf{b}-\mathbf{b}^{\rm tr}\overline{\mathbf{a}})$ as
\begin{align*}
&\big(\mathbf{b}\mathbf{a}\mathbf{b}\overline{\mathbf{a}}\mathbf{b}-(\mathbf{b}^{\rm tr})^2\mathbf{a}\mathbf{b}^{\rm tr}\overline{\mathbf{a}}\big)
+\big(\mathbf{b}\mathbf{b}^{\rm tr}\mathbf{a}\mathbf{b}^{\rm tr}\overline{\mathbf{a}}-\mathbf{b}\mathbf{a}\mathbf{b}\mathbf{b}^{\rm tr}\overline{\mathbf{a}}\big) \\
&+\big(\mathbf{b}^{\rm tr}\mathbf{a}\mathbf{b}\mathbf{b}^{\rm tr}\overline{\mathbf{a}}-\mathbf{b}^{\rm tr}\mathbf{a}\mathbf{b}\overline{\mathbf{a}}\mathbf{b}\big)
+\big((\mathbf{b}^{\rm tr})^2\mathbf{b}-\mathbf{b}\mathbf{b}^{\rm tr}\mathbf{b}\big),
\end{align*}
and using ${\rm tr}(\mathbf{z}^{\rm tr})={\rm tr}(\mathbf{z})$, we obtain
\begin{align*}
{\rm tr}(\mathbf{q})
={\rm tr}(\mathbf{b}\mathbf{a}\mathbf{b}\overline{\mathbf{a}}\mathbf{b})-{\rm tr}(\overline{\mathbf{a}}\mathbf{b}\mathbf{a}\mathbf{b}^2)=t_{12\overline{1}\overline{2}}-t_{21\overline{2}\overline{1}}.
\end{align*}
Then the result follows from (\ref{eq:det=tr}).
\end{proof}

\begin{defn}
\rm Call a pair $(\mathbf{a},\mathbf{b})$ {\it irreducible} if $\sigma_{\mathbf{a},\mathbf{b}}:F_2\to{\rm SL}(3,\mathbb{C})$ is irreducible, or equivalently, neither $\mathbf{a},\mathbf{b}$ nor $\mathbf{a}^{\rm tr},\mathbf{b}^{\rm tr}$ share an eigenvector. Otherwise, call $(\mathbf{a},\mathbf{b})$ {\it reducible}.
\end{defn}

\begin{lem}\label{lem:irr-ordinary}
If $(\mathbf{a},\mathbf{b})$ is irreducible, then
\begin{enumerate}
  \item[\rm(i)] at least one of $\mathbf{a},\mathbf{b}$ is ordinary;
  \item[\rm(ii)] any $\mathbf{u}\in\mathcal{M}$ commuting with both $\mathbf{a}$ and $\mathbf{b}$ is scalar, so that ${\rm Cen}(\mathbf{a})\cap{\rm Cen}(\mathbf{b})=\{\mu\mathbf{e}\colon\mu\in\langle\omega\rangle\}$.
\end{enumerate}
\end{lem}
\begin{proof}
By Burnside's Theorem (see \cite{Bu05} p. 433), ${\rm Im}(\sigma_{\mathbf{a},\mathbf{b}})$ spans $\mathcal{M}$.

(i) Assume $\mathbf{a},\mathbf{b}$ are both special. Then from Corollary \ref{cor:one-is-special} and (\ref{eq:identity2}) we see
${\rm Im}(\sigma_{\mathbf{a},\mathbf{b}})\subset\mathbb{C}\langle\mathbf{e},\mathbf{a},\mathbf{b},\mathbf{a}\mathbf{b},
\mathbf{b}\mathbf{a}\rangle$, so ${\rm Im}(\sigma_{\mathbf{a},\mathbf{b}})$ could not span $\mathcal{M}$.

(ii) If $\mathbf{u}\mathbf{a}=\mathbf{a}\mathbf{u}$ and $\mathbf{u}\mathbf{b}=\mathbf{b}\mathbf{u}$, then $\mathbf{u}\mathbf{w}=\mathbf{w}\mathbf{u}$ for any word $\mathbf{w}$ in $\mathbf{a},\mathbf{b}$; by (i), $\mathbf{u}$ is central in $\mathcal{M}$.
\end{proof}

\begin{prop}
The following conditions are equivalent to each other:
{\rm(i)} $t_{12\overline{12}}=t_{21\overline{21}}$;
{\rm(ii)} 1 is an eigenvalue of $[\mathbf{a},\mathbf{b}]$;
{\rm(iii)} $\det(\mathbf{a}\mathbf{b}-\mathbf{b}\mathbf{a})=0$.
\end{prop}
\begin{proof}
Notice that, for any $\mathbf{z}\in{\rm SL}(3,\mathbb{C})$, ${\rm tr}(\mathbf{z})={\rm tr}(\overline{\mathbf{z}})$ if and only if $1$ is an eigenvalue of $\mathbf{z}$, which is equivalent to the existence of a vector $\mathfrak{a}\ne 0$ such that $\mathbf{z}\mathfrak{a}=\mathfrak{a}$. For $\mathbf{z}=[\mathbf{a},\mathbf{b}]$, this means $(\mathbf{a}\mathbf{b}-\mathbf{b}\mathbf{a})(\overline{\mathbf{a}}\overline{\mathbf{b}}\mathfrak{a})=0$, and the existence of such $\mathfrak{a}$ is equivalent to $\det(\mathbf{a}\mathbf{b}-\mathbf{b}\mathbf{a})=0$.
\end{proof}

\begin{rmk}
\rm Interestingly, as shown in \cite{PauWi17} Theorem 4.1 and 4.2, when $\mathbf{a},\mathbf{b}\in{\rm SU}(2,1)$, the property that 1 is an eigenvalue of $[\mathbf{a},\mathbf{b}]$ gives geometric information about the action of the subgroup $\langle\mathbf{a},\mathbf{b}\rangle$.
\end{rmk}

\subsection{The asymmetric part}\label{sec:assymetric}

Suppose $t_{12\overline{12}}\ne t_{21\overline{21}}$, then obviously $(\mathbf{a},\mathbf{b})$ is irreducible.
Fix $(\mathbf{a},\mathbf{b})$, and define a bilinear pairing on $\mathcal{M}$ by
$$((\mathbf{u},\mathbf{v}))=(t_{12\overline{12}}-t_{21\overline{21}})^{-1}\cdot{\rm tr}(\mathbf{u}\mathbf{v}).$$
Direct computations lead to the following table (with zeros in empty cases):
\begin{center}
\begin{tabular}{|c|c|c|c|c|c|c|c|c|c|}
  \hline
  $((\mathbf{u},\mathbf{v}))$ & $\mathbf{e}$ & $\mathbf{a}$ & $\overline{\mathbf{a}}$ & $\mathbf{b}$ & $\overline{\mathbf{b}}$ & $\theta(\mathbf{a}\mathbf{b})$ & $\theta(\mathbf{a}\overline{\mathbf{b}})$ & $\theta(\overline{\mathbf{a}}\mathbf{b})$ & $\theta(\overline{\mathbf{a}}\overline{\mathbf{b}})$ \\
  \hline
  $\delta([\mathbf{a},\mathbf{b}])$ & $1$ & $t_1$ & \ & \ & $t_{\overline{2}}$ & $t_{12}$ & $t_{1\overline{2}}+t_1t_{\overline{2}}$ & \ & $t_{\overline{12}}$ \\
  \hline
  $\delta(\mathbf{b}\lrcorner\overline{\mathbf{a}})$ & \ & $1$ & \ & \ & \ & $t_2$ & $t_{\overline{2}}$ & \ & \ \\
  \hline
  $\delta(\mathbf{b}\lrcorner\mathbf{a})$ & \ & \ & $-1$ & \ & \ & \ & \ & $-t_2$ & $-t_{\overline{2}}$ \\
  \hline
  $\delta(\mathbf{a}\lrcorner\overline{\mathbf{b}})$ & \ & \ & \ & $-1$ & \ & $-t_1$ & \ & $-t_{\overline{1}}$ & \ \\
  \hline
  $\delta(\mathbf{a}\lrcorner\mathbf{b})$ & \ & \ & \ & \ & $1$ & \ & $t_1$ & \ & $t_{\overline{1}}$ \\
  \hline
  $\delta(\overline{\mathbf{a}}\overline{\mathbf{b}})$ & \ & \ & \ & \ & \ & $1$ & \ & \ & \ \\
  \hline
  $\delta(\overline{\mathbf{a}}\mathbf{b})$ & \ & \ & \ & \ & \ & \ & $-1$ & \ & \ \\
  \hline
  $\delta(\mathbf{a}\overline{\mathbf{b}})$ & \ & \ & \ & \ & \ & \ & \ & $-1$ & \ \\
  \hline
  $\delta(\mathbf{a}\mathbf{b})$ & \ & \ & \ & \ & \ & \ & \ & \ & $1$ \\
  \hline
\end{tabular} 
\end{center}

\begin{thm}\label{thm:basis-1}
If $t_{12\overline{12}}\ne t_{21\overline{21}}$, then $\mathfrak{B}_\delta$, $\mathfrak{B}_\theta$ are both bases of $\mathcal{M}$, where
\begin{align*}
\mathfrak{B}_\delta&=\big\{\delta(\mathbf{a}\mathbf{b}), \delta(\mathbf{a}\overline{\mathbf{b}}), \delta(\overline{\mathbf{a}}\mathbf{b}),
\delta(\overline{\mathbf{a}}\overline{\mathbf{b}}), \delta(\mathbf{a}\lrcorner\mathbf{b}), \delta(\mathbf{a}\lrcorner\overline{\mathbf{b}}),
\delta(\mathbf{b}\lrcorner\mathbf{a}), \delta(\mathbf{b}\lrcorner\overline{\mathbf{a}}), \delta([\mathbf{a},\mathbf{b}])\big\}, \\
\mathfrak{B}_\theta&=\big\{\mathbf{e}, \mathbf{a}, \overline{\mathbf{a}}, \mathbf{b}, \overline{\mathbf{b}}, \theta(\mathbf{a}\mathbf{b}), \theta(\mathbf{a}\overline{\mathbf{b}}), \theta(\overline{\mathbf{a}}\mathbf{b}), \theta(\overline{\mathbf{a}}\overline{\mathbf{b}})\big\}.
\end{align*}
\end{thm}

\begin{proof}
It is a special case of a general result: if $\langle\cdot,\cdot\rangle$ is a bilinear pairing on a vector space $\mathcal{V}$ of dimension $n$, and $\mathfrak{B}_1=\{\mathfrak{u}_1,\ldots,\mathfrak{u}_n\}, \mathfrak{B}_2=\{\mathfrak{v}_1,\ldots,\mathfrak{v}_n\}$ are subsets of $\mathcal{V}$ such that the matrix $P=[\langle \mathfrak{u}_i,\mathfrak{v}_j\rangle]_{n\times n}$ is invertible, then $\mathfrak{B}_1,\mathfrak{B}_2$ are both bases of $\mathcal{V}$.

To show this, suppose $\sum_{i=1}^na_i\mathfrak{u}_i=0$, with $a_i\in\mathbb{C}$. Then $\sum_{i=1}^na_i\langle\mathfrak{u}_i,\mathfrak{v}_j\rangle=0$ for each $j$. Hence $[a_1,\ldots,a_n]P=0$, implying $a_1=\ldots=a_n=0$. Similarly, any linear relation $\sum_{j=1}^nb_j\mathfrak{v}_j=0$ holds only when $b_1=\cdots=b_n=0$.
\end{proof}

\begin{rmk}\label{rmk:commutators}
\rm In $\mathfrak{B}_\delta$, one may replace $\delta([\mathbf{a},\mathbf{b}])$ by $\delta([\mathbf{b},\mathbf{a}])$. Indeed,
\begin{align*}
\delta([\mathbf{a},\mathbf{b}])+\delta([\mathbf{b},\mathbf{a}])=\ &
(t_{\overline{1}2}-t_{\overline{1}}t_2)\delta(\mathbf{a}\overline{\mathbf{b}})+(t_1t_{\overline{2}}-t_{1\overline{2}})\delta(\overline{\mathbf{a}}\mathbf{b}) \\
&+t_{\overline{2}}\delta(\mathbf{a}\lrcorner\mathbf{b})+t_2\delta(\mathbf{a}\lrcorner\overline{\mathbf{b}})+t_{\overline{1}}\delta(\mathbf{b}\lrcorner\mathbf{a})+t_1\delta(\mathbf{b}\lrcorner\overline{\mathbf{a}}).
\end{align*}
This is obtained by writing
$$\delta([\mathbf{a},\mathbf{b}])+\delta([\mathbf{b},\mathbf{a}])=\theta(\mathbf{a}\lrcorner\mathbf{b})\overline{\mathbf{b}}-\overline{\mathbf{b}}\theta(\mathbf{a}\lrcorner\mathbf{b})
+\theta(\mathbf{b}\lrcorner\mathbf{a})\overline{\mathbf{a}}-\overline{\mathbf{a}}\theta(\mathbf{b}\lrcorner\mathbf{a})$$
and applying (\ref{eq:identity1}).
\end{rmk}

\subsection{The symmetric slice}  \label{sec:symmetric}

\begin{thm} \label{thm:symmetric}
For any irreducible pair $(\mathbf{a},\mathbf{b})$, there exists $\mathbf{c}\in{\rm SL}(3,\mathbb{C})$ with $\mathbf{c}\lrcorner\mathbf{a}$, $\mathbf{c}\lrcorner\mathbf{b}$ both symmetric if and only if $t_{12\overline{12}}=t_{21\overline{21}}$.
\end{thm}

\begin{proof}
The ``only if" part is trivial: if $\mathbf{c}\lrcorner\mathbf{a}, \mathbf{c}\lrcorner\mathbf{b}\in\mathcal{M}^0$, then
$$t_{12\overline{12}}={\rm tr}([\mathbf{c}\lrcorner\mathbf{a},\mathbf{c}\lrcorner\mathbf{b}])={\rm tr}([\mathbf{c}\lrcorner\mathbf{a},\mathbf{c}\lrcorner\mathbf{b}]^{\rm tr})=t_{21\overline{21}}.$$

For the ``if" part, suppose $t_{12\overline{12}}=t_{21\overline{21}}$.
By Lemma \ref{lem:irr-ordinary} and the symmetry between $\mathbf{a},\mathbf{b}$, we can assume that $\mathbf{a}$ is ordinary; also we can assume $\mathbf{a}$ to be symmetric.
By Lemma \ref{lem:tr-tr}, there exists a nontrivial linear relation
$$\alpha(\mathbf{b}-\mathbf{b}^{\rm tr})^\neg+\beta(\mathbf{a}\mathbf{b}-\mathbf{b}^{\rm tr}\mathbf{a})^\neg+\gamma(\overline{\mathbf{a}}\mathbf{b}-\mathbf{b}^{\rm tr}\overline{\mathbf{a}})^\neg=\mathbf{0};$$
equivalently, $\tilde{\mathbf{c}}\mathbf{b}=\mathbf{b}^{\rm tr}\tilde{\mathbf{c}}$, where $\tilde{\mathbf{c}}=\alpha\mathbf{e}+\beta\mathbf{a}+\gamma\overline{\mathbf{a}}\in\mathcal{M}$.

We show $\det(\tilde{\mathbf{c}})\ne 0$.
Assume on the contrary that $\det(\tilde{\mathbf{c}})=0$. Suppose
$\det(t\mathbf{e}-\mathbf{a})=(t-\lambda_1)(t-\lambda_2)(t-\lambda_3).$
Let $\mathcal{E}(\lambda_i)$ denote the eigenspace of $\mathbf{a}$ with respect to $\lambda_i$, which is one-dimensional.
Let $f(t)=\beta t^2+\alpha t+\gamma$. Then $f(\lambda_i)=0$ for at least one $i$, and $f(\mathbf{a})=\mathbf{a}\tilde{\mathbf{c}}$.
\begin{itemize}
  \item If the greatest common divisor of $f(t)$ and $\det(t\mathbf{e}-\mathbf{a})$ is $t-\lambda_i$ for some $i$, then taking a nonzero $\mathfrak{a}\in \mathcal{E}(\lambda_i)$, we have
      $$f(\mathbf{a})(\mathbf{b}\mathfrak{a})=\mathbf{a}\tilde{\mathbf{c}}\mathbf{b}\mathfrak{a}=\mathbf{a}\mathbf{b}^{\rm tr}\tilde{\mathbf{c}}\mathfrak{a}=\mathbf{a}\mathbf{b}^{\rm tr}\overline{\mathbf{a}}f(\mathbf{a})\mathfrak{a}=0.$$
      Consequently, $(\mathbf{a}-\lambda_i\mathbf{e})(\mathbf{b}\mathfrak{a})=0$, i.e. $\mathbf{b}\mathfrak{a}\in \mathcal{E}(\lambda_i)$, showing that $\mathbf{b}\mathfrak{a}$ is a multiple of $\mathfrak{a}$, so $\mathfrak{a}$ is a common eigenvector of $\mathbf{a}$ and $\mathbf{b}$. This contradicts the irreducibility of $(\mathbf{a},\mathbf{b})$.
  \item If $f(t)=\beta(t-\lambda_i)(t-\lambda_j)$ for $i\ne j$, (note that $\lambda_i=\lambda_j$ is allowed), then $\{\tilde{\mathbf{c}}\mathfrak{c}\colon\mathfrak{c}\in\mathbb{C}^3\}=\mathcal{E}(\lambda_k)$, where $\{k\}=\{1,2,3\}-\{i,j\}$. Taking $\mathfrak{b}$ with $\mathfrak{a}:=\tilde{\mathbf{c}}\mathfrak{b}\ne0$,
      we have
      $$\mathbf{b}^{\rm tr}\mathfrak{a}=\mathbf{b}^{\rm tr}(\tilde{\mathbf{c}}\mathfrak{b})=\tilde{\mathbf{c}}(\mathbf{b}\mathfrak{b})\in \mathcal{E}(\lambda_k)=\mathbb{C}\langle\mathfrak{a}\rangle,$$
      so $\mathfrak{a}$ is an eigenvector of both $\mathbf{a}$ and $\mathbf{b}^{\rm tr}$. This also contradicts the irreducibility of $(\mathbf{a},\mathbf{b})$.
\end{itemize}

Thus, $\det(\tilde{\mathbf{c}})\ne 0$.
Taking $\kappa$ with $\kappa^2=\det(\tilde{\mathbf{c}})$, we have $\kappa^{-1}\tilde{\mathbf{c}}\in{\rm Cen}(\mathbf{a})$. By Lemma \ref{lem:centralizer-ordinary}, there exists $\mathbf{c}\in{\rm Cen}(\mathbf{a})\subset\mathbb{C}\langle\mathbf{e},\mathbf{a},\overline{\mathbf{a}}\rangle\subset\mathcal{M}^0$ with $\mathbf{c}^2=\kappa^{-1}\tilde{\mathbf{c}}$, which ensures $\mathbf{c}\mathbf{b}\overline{\mathbf{c}}=\overline{\mathbf{c}}\mathbf{b}^{\rm tr}\mathbf{c}$, i.e. $\mathbf{c}\lrcorner\mathbf{b}\in\mathcal{M}^0$.
\end{proof}

Let $\mathcal{R}^{\rm irr,sym}_{{\rm SL}(3,\mathbb{C})}(F_2)$ denote the set of irreducible pairs $(\mathbf{a},\mathbf{b})$ with $\mathbf{a},\mathbf{b}\in\mathcal{M}^0$.
Let $\mathcal{X}^{\rm irr,sym}_{{\rm SL}(3,\mathbb{C})}(F_2)=\big\{\chi_\sigma\colon\sigma\in\mathcal{R}^{\rm irr,sym}_{{\rm SL}(3,\mathbb{C})}(F_2)\big\}$, which coincides with the subset of $\mathcal{X}_{{\rm SL}(3,\mathbb{C})}(F_2)$ consisting of irreducible characters with $t_{12\overline{12}}=t_{21\overline{21}}$.

\begin{thm}
The action of ${\rm SO}(3,\mathbb{C})$ on $\mathcal{R}^{\rm irr,sym}_{{\rm SL}(3,\mathbb{C})}(F_2)$ via conjugation is free, and is transitive on each fiber of $\chi:\mathcal{R}^{\rm irr,sym}_{{\rm SL}(3,\mathbb{C})}(F_2)\to\mathcal{X}^{\rm irr,sym}_{{\rm SL}(3,\mathbb{C})}(F_2)$.
Consequently,
$$\mathcal{X}^{\rm irr,sym}_{{\rm SL}(3,\mathbb{C})}(F_2)
\cong\mathcal{R}^{\rm irr,sym}_{{\rm SL}(3,\mathbb{C})}(F_2)/{\rm SO}(3,\mathbb{C}).$$
\end{thm}

\begin{proof}
The first assertion follows from Lemma \ref{lem:irr-ordinary} (ii) and the fact that ${\rm SO}(3,\mathbb{C})\cap\langle\omega\rangle\mathbf{e}=\{\mathbf{e}\}$.

For the second one, suppose $(\mathbf{a},\mathbf{b}),(\mathbf{a}',\mathbf{b}')\in\mathcal{R}^{\rm irr,sym}_{{\rm SL}(3,\mathbb{C})}(F_2)$ have the same character. Then there exists $\mathbf{c}\in{\rm SL}(3,\mathbb{C})$ such that $\mathbf{c}\lrcorner\mathbf{a}=\mathbf{a}'$, $\mathbf{c}\lrcorner\mathbf{b}=\mathbf{b}'$. From $\mathbf{c}\lrcorner\mathbf{a}=\mathbf{a}'=(\mathbf{a}')^{\rm tr}=(\mathbf{c}\lrcorner\mathbf{a})^{\rm tr}$ we deduce $\mathbf{c}^{\rm tr}\mathbf{c}\in{\rm Cen}(\mathbf{a})$, and from $\mathbf{c}\lrcorner\mathbf{b}=(\mathbf{c}\lrcorner\mathbf{b})^{\rm tr}$ we deduce $\mathbf{c}^{\rm tr}\mathbf{c}\in{\rm Cen}(\mathbf{b})$. Hence by Lemma \ref{lem:irr-ordinary} (ii), $\mathbf{c}^{\rm tr}\mathbf{c}=\mu\mathbf{e}$ for some $\mu\in\langle\omega\rangle$. Then $\mu\mathbf{c}\in{\rm SO}(3,\mathbb{C})$, and $(\mu\mathbf{c})\lrcorner\mathbf{a}=\mathbf{a}'$, $(\mu\mathbf{c})\lrcorner\mathbf{b}=\mathbf{b}'$.
\end{proof}

\medskip

We aim to find convenient bases for $\mathcal{M}^0$, $\mathcal{M}^1$, $\mathcal{M}$.

\begin{lem}\label{lem:span}
If $\mathbf{a},\mathbf{b}$ are symmetric and the pair $(\mathbf{a},\mathbf{b})$ is irreducible, then
\begin{align*}
\mathcal{M}^0&=\mathbb{C}\langle\mathbf{e},\mathbf{a},\overline{\mathbf{a}},\mathbf{b},\overline{\mathbf{b}},\theta(\mathbf{a}\mathbf{b}),\theta(\mathbf{a}\overline{\mathbf{b}}),\theta(\overline{\mathbf{a}}\mathbf{b}),
\theta(\overline{\mathbf{a}}\overline{\mathbf{b}}),\theta([\mathbf{a},\mathbf{b}]),\theta([\mathbf{b},\mathbf{a}])\rangle, \\
\mathcal{M}^1&=\mathbb{C}\langle\delta(\mathbf{a}\mathbf{b}),\delta(\overline{\mathbf{a}}\mathbf{b}),\delta(\mathbf{a}\overline{\mathbf{b}}),
\delta(\overline{\mathbf{a}}\overline{\mathbf{b}}),
\delta(\mathbf{a}\lrcorner\mathbf{b}),\delta(\mathbf{a}\lrcorner\overline{\mathbf{b}}),\delta(\mathbf{b}\lrcorner\mathbf{a}),\delta(\mathbf{b}\lrcorner\overline{\mathbf{a}}),\delta([\mathbf{a},\mathbf{b}])\rangle.
\end{align*}
\end{lem}
\begin{proof}
Note the basic fact: 
if $\mathcal{M}$ is spanned by $\mathbf{x}_1,\ldots,\mathbf{x}_n$, then
$\mathcal{M}^0$ is spanned by $\mathbf{x}_1+\mathbf{x}_1^{\rm tr},\ldots,\mathbf{x}_n+\mathbf{x}_n^{\rm tr}$, and $\mathcal{M}^1$ is spanned by $\mathbf{x}_1-\mathbf{x}_1^{\rm tr},\ldots,\mathbf{x}_n-\mathbf{x}_n^{\rm tr}$.

The first equation follows from Proposition \ref{prop:express} together with (\ref{eq:identity1}), and the second one follows from Proposition \ref{prop:express} together with Remark \ref{rmk:commutators}.
\end{proof}

\begin{thm} \label{thm:basis-2}
Suppose $\mathbf{a},\mathbf{b}$ are symmetric and the pair $(\mathbf{a},\mathbf{b})$ is irreducible.

{\rm(i)} If $\mathbf{a}$ is special (so that $\mathbf{b}$ is ordinary), then $\{\mathbf{e},\mathbf{a},\mathbf{b},\overline{\mathbf{b}},\theta(\mathbf{a}\mathbf{b}),\theta(\mathbf{a}\overline{\mathbf{b}})\}$ is a basis of $\mathcal{M}^0$,
and $\{\delta(\mathbf{a}\mathbf{b}), \delta(\mathbf{a}\overline{\mathbf{b}}),\delta(\mathbf{b}\lrcorner\mathbf{a})\}$ is a basis of $\mathcal{M}^1$.

{\rm(ii)} If $\mathbf{a},\mathbf{b}$ are both ordinary, then at least one of the following occurs:
\begin{enumerate}
  \item[\rm 1)] $\{\mathbf{e},\mathbf{a},\overline{\mathbf{a}},\mathbf{b},\overline{\mathbf{b}},\theta(\mathbf{a}\mathbf{b})\}$ is a basis of $\mathcal{M}^0$,
      and $\{\delta(\mathbf{a}\mathbf{b}), \delta(\overline{\mathbf{a}}\mathbf{b}), \delta(\mathbf{a}\overline{\mathbf{b}})\}$ is a basis of $\mathcal{M}^1$,
      so $\mathfrak{B}^{\rm sym}_1=\{\mathbf{e},\mathbf{a},\overline{\mathbf{a}},\mathbf{b},\overline{\mathbf{b}},
      \mathbf{a}\mathbf{b}, \mathbf{b}\mathbf{a}, \overline{\mathbf{a}}\mathbf{b},\mathbf{a}\overline{\mathbf{b}}\}$ is a basis of $\mathcal{M}$;
  \item[\rm 2)] $\{\mathbf{e},\mathbf{a},\overline{\mathbf{a}},\mathbf{b},\overline{\mathbf{b}},
      \theta(\overline{\mathbf{a}}\mathbf{b})\}$ is a basis of $\mathcal{M}^0$,
      and $\{\delta(\mathbf{a}\mathbf{b}), \delta(\overline{\mathbf{a}}\mathbf{b}), \delta(\overline{\mathbf{a}}\overline{\mathbf{b}})\}$ is a basis of $\mathcal{M}^1$,
      so $\mathfrak{B}^{\rm sym}_2=\{\mathbf{e},\mathbf{a},\overline{\mathbf{a}},\mathbf{b},\overline{\mathbf{b}},
      \mathbf{a}\mathbf{b}, \mathbf{b}\overline{\mathbf{a}}, \overline{\mathbf{a}}\mathbf{b},\overline{\mathbf{a}}\overline{\mathbf{b}}\}$ is a basis of $\mathcal{M}$.
\end{enumerate}
Therefore, in any case,
\begin{align*}
\mathcal{M}^0&=\mathbb{C}\langle\mathbf{e},\mathbf{a},\overline{\mathbf{a}},\mathbf{b},\overline{\mathbf{b}},
\theta(\mathbf{a}\mathbf{b}),\theta(\overline{\mathbf{a}}\mathbf{b})\rangle, \\
\mathcal{M}^1&=\mathbb{C}\langle\delta(\mathbf{a}\mathbf{b}),\delta(\overline{\mathbf{a}}\mathbf{b}),
\delta(\mathbf{a}\overline{\mathbf{b}}),\delta(\overline{\mathbf{a}}\overline{\mathbf{b}})\rangle, \\
\mathcal{M}&=\mathbb{C}\langle\mathbf{e},\mathbf{a},\overline{\mathbf{a}},\mathbf{b},\overline{\mathbf{b}},
\mathbf{a}\mathbf{b},\mathbf{b}\mathbf{a},\overline{\mathbf{a}}\mathbf{b},\mathbf{b}\overline{\mathbf{a}},\mathbf{a}\overline{\mathbf{b}},
\overline{\mathbf{a}}\overline{\mathbf{b}}\rangle.
\end{align*}
\end{thm}

\begin{proof}
(i) This is an immediate consequence of Corollary \ref{cor:one-is-special} and the basic fact stated in the proof of Lemma \ref{lem:span}.

(ii) The result is established by proving a series of assertions.
\begin{enumerate}
  \item $\mathbf{e},\mathbf{a},\overline{\mathbf{a}},\mathbf{b},\overline{\mathbf{b}}$ are linearly independent.

        Suppose
        $$\alpha_0\mathbf{e}+\alpha_1\mathbf{a}+\alpha_2\overline{\mathbf{a}}+\beta_1\mathbf{b}+\beta_2\overline{\mathbf{b}}=\mathbf{0}.$$
        Then $\beta_1\mathbf{b}+\beta_2\overline{\mathbf{b}}$ commutes with $\mathbf{a}$ and $\mathbf{b}$, so by Lemma \ref{lem:irr-ordinary} (ii), $\beta_1\mathbf{b}+\beta_2\overline{\mathbf{b}}=\gamma\mathbf{e}$ for some $\gamma$. Since $\mathbf{b}$ is ordinary, we have $\beta_1=\beta_2=\gamma=0$, so that $\alpha_0\mathbf{e}+\alpha_1\mathbf{a}+\alpha_2\overline{\mathbf{a}}=\mathbf{0}$. This forces $\alpha_0=\alpha_1=\alpha_2=0$, as $\mathbf{a}$ is ordinary.
  \item Let $\mathcal{V}^0=\mathbb{C}\langle\mathbf{e},\mathbf{a},\overline{\mathbf{a}},\mathbf{b},\overline{\mathbf{b}}\rangle$.
        Then either $\theta(\mathbf{a}\mathbf{b})\notin\mathcal{V}^0$ or $\theta(\overline{\mathbf{a}}\mathbf{b})\notin\mathcal{V}^0$.

        Assume on the contrary that $\theta(\mathbf{a}\mathbf{b}),\theta(\overline{\mathbf{a}}\mathbf{b})\in\mathcal{V}^0$.
        Then
        $$\theta(\mathbf{a}^2\mathbf{b})=t_1\theta(\mathbf{a}\mathbf{b})-2t_{\overline{1}}\mathbf{b}+\theta(\overline{\mathbf{a}}\mathbf{b})\in\mathcal{V}^0;$$
        by (\ref{eq:identity1}), $\theta(\mathbf{a}\lrcorner\mathbf{b})\in\mathcal{V}^0$; by (\ref{eq:identity2}), $\mathbf{a}\mathbf{b}\mathbf{a}\in\mathcal{V}^0$.
        Suppose
        \begin{align}
        \theta(\mathbf{a}\mathbf{b})=\alpha_0\mathbf{e}+\alpha_1\mathbf{a}+\alpha_2\overline{\mathbf{a}}
        +\beta_1\mathbf{b}+\beta_2\overline{\mathbf{b}}.  \label{eq:linear-comb-1}
        \end{align}

        Multiplying this by $\mathbf{b}$ on the left and right, we respectively obtain
        \begin{align}
        \mathbf{b}\mathbf{a}\mathbf{b}+\mathbf{b}^2\mathbf{a}=\alpha_0\mathbf{b}+\alpha_1\mathbf{b}\mathbf{a}+\alpha_2\mathbf{b}\overline{\mathbf{a}}+\beta_1\mathbf{b}^2+\beta_2\mathbf{e},  \label{eq:consequence1} \\
        \mathbf{a}\mathbf{b}^2+\mathbf{b}\mathbf{a}\mathbf{b}=\alpha_0\mathbf{b}+\alpha_1\mathbf{a}\mathbf{b}+\alpha_2\overline{\mathbf{a}}\mathbf{b}+\beta_1\mathbf{b}^2+\beta_2\mathbf{e}.  \label{eq:consequence2}
        \end{align}
        Their sum yields $\theta(\mathbf{a}\mathbf{b}^2)+2\mathbf{b}\mathbf{a}\mathbf{b}\in\mathcal{V}^0$. By (\ref{eq:identity2}), $\theta(\mathbf{a}\mathbf{b}^2)\in\mathcal{V}^0$. Furthermore, since $\theta(\mathbf{a}\mathbf{b}^2)=t_2\theta(\mathbf{a}\mathbf{b})-2t_{\overline{2}}\mathbf{a}+\theta(\mathbf{a}\overline{\mathbf{b}})$, we have $\theta(\mathbf{a}\overline{\mathbf{b}})\in\mathcal{V}^0$.

        A similar argument leads to $\theta(\overline{\mathbf{a}}\mathbf{b}^2),\theta(\overline{\mathbf{a}}\overline{\mathbf{b}})\in\mathcal{V}^0$, and moreover, $$\theta(\overline{\mathbf{a}}^2\overline{\mathbf{b}})=t_{\overline{1}}\theta(\overline{\mathbf{a}}\overline{\mathbf{b}})
        -2t_{1}\overline{\mathbf{b}}+\theta(\mathbf{a}\overline{\mathbf{b}})\in\mathcal{V}^0.$$
        Also by (\ref{eq:identity1}), $\theta(\mathbf{a}\lrcorner\overline{\mathbf{b}}),\theta(\mathbf{b}\lrcorner\mathbf{a}),\theta(\mathbf{b}\lrcorner\overline{\mathbf{a}})\in\mathcal{V}^0$.

        Multiplying (\ref{eq:linear-comb-1}) by $\overline{\mathbf{a}}\overline{\mathbf{b}}$ on the right and by $\overline{\mathbf{b}}\overline{\mathbf{a}}$ on the left, and then adding the resulted equations, we obtain
        $$\theta([\mathbf{a},\mathbf{b}])+2\mathbf{e}=\alpha_0\theta(\overline{\mathbf{a}}\overline{\mathbf{b}})+2\alpha_1\overline{\mathbf{b}}+\alpha_2\theta(\overline{\mathbf{a}}^2\overline{\mathbf{b}})
        +\beta_1\theta(\mathbf{b}\lrcorner\overline{\mathbf{a}})+2\beta_2\overline{\mathbf{b}}\overline{\mathbf{a}}\overline{\mathbf{b}},$$
        implying $\theta([\mathbf{a},\mathbf{b}])\in\mathcal{V}^0$. Similarly, $\theta([\mathbf{b},\mathbf{a}])\in\mathcal{V}^0$.

        By Lemma \ref{lem:span}, $\mathcal{M}^0=\mathcal{V}^0$ which is absurd, so the assumption $\theta(\mathbf{a}\mathbf{b}),\theta(\overline{\mathbf{a}}\mathbf{b})\in\mathcal{V}^0$ is false.

        Therefore,
        $\mathcal{M}^0=\mathbb{C}\langle\mathbf{e},\mathbf{a},\overline{\mathbf{a}},\mathbf{b},\overline{\mathbf{b}},\theta(\mathbf{a}\mathbf{b})\rangle$ or $\mathcal{M}^0=\mathbb{C}\langle\mathbf{e},\mathbf{a},\overline{\mathbf{a}},\mathbf{b},\overline{\mathbf{b}},
        \theta(\overline{\mathbf{a}}\mathbf{b})\rangle$.
  \item Let $\mathcal{V}^1=\mathbb{C}\langle\delta(\overline{\mathbf{a}}\mathbf{b}),\delta(\mathbf{a}\overline{\mathbf{b}})\rangle$. Then either $\delta(\mathbf{a}\mathbf{b})\notin\mathcal{V}^1$ or $\delta(\overline{\mathbf{a}}\overline{\mathbf{b}})\notin\mathcal{V}^1$.

        Assume on the contrary that $\delta(\mathbf{a}\mathbf{b}),\delta(\overline{\mathbf{a}}\overline{\mathbf{b}})\in\mathcal{V}^1$.
        Suppose
        \begin{align}
        \delta(\mathbf{a}\mathbf{b})=\eta_1\delta(\overline{\mathbf{a}}\mathbf{b})+\eta_2\delta(\mathbf{a}\overline{\mathbf{b}}).  \label{eq:linear-comb-2}
        \end{align}
        Note that $\eta_1\ne 0$: otherwise, $\mathbf{a}(\mathbf{b}-\eta_2\overline{\mathbf{b}})=(\mathbf{b}-\eta_2\overline{\mathbf{b}})\mathbf{a}$, so that $\mathbf{b}-\eta_2\overline{\mathbf{b}}\in{\rm Cen}(\mathbf{a})\cap{\rm Cen}(\mathbf{b})$, which would contradict the assumption that $\mathbf{b}$ is ordinary. Similarly, $\eta_2\ne 0$.

        Multiplying (\ref{eq:linear-comb-2}) by $\mathbf{a}$ on the left and right, and then taking sum, we obtain $\delta(\mathbf{a}^2\mathbf{b})=-\eta_1\delta(\mathbf{a}\lrcorner\mathbf{b})+\eta_2\delta(\mathbf{a}^2\overline{\mathbf{b}})$; since $\delta(\mathbf{a}^2\mathbf{b}^{\pm1})=t_1\delta(\mathbf{a}\mathbf{b}^{\pm1})+\delta(\overline{\mathbf{a}}\mathbf{b}^{\pm1})$, we have $\delta(\mathbf{a}\lrcorner\mathbf{b})\in\mathcal{V}^1$.

        Multiplying (\ref{eq:linear-comb-2}) by $\overline{\mathbf{a}}$ on the left and right, and then taking sum, we obtain $\delta(\mathbf{a}\lrcorner\mathbf{b})=\eta_1\delta(\overline{\mathbf{a}}^2\mathbf{b})+\eta_2\delta(\mathbf{a}\lrcorner\overline{\mathbf{b}})$, implying $\delta(\mathbf{a}\lrcorner\overline{\mathbf{b}})\in\mathcal{V}^1$.

        Similarly, $\delta(\mathbf{b}\lrcorner\mathbf{a}),\delta(\mathbf{b}\lrcorner\overline{\mathbf{a}})\in\mathcal{V}^1$.

        Multiplying (\ref{eq:linear-comb-2}) by $\overline{\mathbf{a}}\overline{\mathbf{b}}$ on the right and $\overline{\mathbf{b}}\overline{\mathbf{a}}$ on the left, and then taking sum, we obtain $\delta([\mathbf{a},\mathbf{b}])=\eta_1\delta(\overline{\mathbf{a}}\mathbf{b}\overline{\mathbf{a}}\overline{\mathbf{b}})+\eta_2\delta(\mathbf{a}\overline{\mathbf{b}}\overline{\mathbf{a}}\overline{\mathbf{b}})$.
        Since
        \begin{align*}
        \delta(\overline{\mathbf{a}}\mathbf{b}\overline{\mathbf{a}}\overline{\mathbf{b}})&\stackrel{(\ref{eq:identity3})}
        =\delta\big((t_{\overline{1}2}\overline{\mathbf{a}}-t_{1\overline{2}}\overline{\mathbf{b}}
        +\overline{\mathbf{b}}\mathbf{a}\overline{\mathbf{b}})\overline{\mathbf{b}}\big)=t_{\overline{1}2}\delta(\overline{\mathbf{a}}\overline{\mathbf{b}})+\delta(\overline{\mathbf{b}}\mathbf{a}\overline{\mathbf{b}}^2) \\
        &=t_{\overline{1}2}\delta(\overline{\mathbf{a}}\overline{\mathbf{b}})+\delta\big(\overline{\mathbf{b}}\mathbf{a}(t_{\overline{2}}\overline{\mathbf{b}}-t_2\mathbf{e}+\mathbf{b})\big) \\
        &=t_{\overline{1}2}\delta(\overline{\mathbf{a}}\overline{\mathbf{b}})+t_2\delta(\mathbf{a}\overline{\mathbf{b}})-\delta(\mathbf{b}\lrcorner\mathbf{a})\in\mathcal{V}^1, \\
        \delta(\mathbf{a}\overline{\mathbf{b}}\overline{\mathbf{a}}\overline{\mathbf{b}})&\stackrel{(\ref{eq:identity3})}
        =\delta\big(\mathbf{a}(t_{\overline{12}}\overline{\mathbf{b}}-t_{12}\mathbf{a}+\mathbf{a}\mathbf{b}\mathbf{a})\big)
        =t_{\overline{12}}\delta(\mathbf{a}\overline{\mathbf{b}})+\delta(\mathbf{a}^2\mathbf{b}\mathbf{a}) \\
        &=t_{\overline{12}}\delta(\mathbf{a}\overline{\mathbf{b}})+\delta((t_1\mathbf{a}-t_{\overline{1}}\mathbf{e}+\overline{\mathbf{a}})\mathbf{b}\mathbf{a}\big)  \\
        &=t_{\overline{12}}\delta(\mathbf{a}\overline{\mathbf{b}})+t_{\overline{1}}\delta(\mathbf{a}\mathbf{b})-\delta(\mathbf{a}\lrcorner\mathbf{b})\in\mathcal{V}^1,
        \end{align*}
        we have $\delta([\mathbf{a},\mathbf{b}])\in\mathcal{V}^1$.

        By Lemma \ref{lem:span}, $\mathcal{M}^1=\mathcal{V}^1$ which is absurd.

        Therefore,
        $\mathcal{M}^1=\mathbb{C}\langle\delta(\mathbf{a}\mathbf{b}),\delta(\overline{\mathbf{a}}\mathbf{b}),
        \delta(\mathbf{a}\overline{\mathbf{b}})\rangle$ or        $\mathcal{M}^1=\mathbb{C}\langle\delta(\overline{\mathbf{a}}\mathbf{b}),\delta(\mathbf{a}\overline{\mathbf{b}}),
        \delta(\overline{\mathbf{a}}\overline{\mathbf{b}})\rangle$.
  \item Finally, if $\mathcal{M}^1=\mathbb{C}\langle\delta(\mathbf{a}\mathbf{b}),\delta(\overline{\mathbf{a}}\mathbf{b}),\delta(\mathbf{a}\overline{\mathbf{b}})\rangle$, then
        $\mathcal{M}^0=\mathbb{C}\langle\mathbf{e},\mathbf{a},\overline{\mathbf{a}},\mathbf{b},\overline{\mathbf{b}},\theta(\mathbf{a}\mathbf{b})\rangle$: otherwise, supposing (\ref{eq:linear-comb-1}), the difference between (\ref{eq:consequence1}) and (\ref{eq:consequence2}) would imply
        $$\alpha_1\delta(\mathbf{a}\mathbf{b})+\alpha_2\delta(\overline{\mathbf{a}}\mathbf{b})=\delta(\mathbf{a}\mathbf{b}^2)
        =t_1\delta(\mathbf{a}\mathbf{b})+\delta(\mathbf{a}\overline{\mathbf{b}}),$$
        so that $\dim\mathcal{M}^1\le 2$ which is absurd.

        Similarly, if $\mathcal{M}^1=\mathbb{C}\langle\delta(\overline{\mathbf{a}}\mathbf{b}),\delta(\mathbf{a}\overline{\mathbf{b}}),\delta(\overline{\mathbf{a}}\overline{\mathbf{b}})\rangle$, then
        $\mathcal{M}^0=\mathbb{C}\langle\mathbf{e},\mathbf{a},\overline{\mathbf{a}},\mathbf{b},\overline{\mathbf{b}},\theta(\overline{\mathbf{a}}\mathbf{b})\rangle$.
\end{enumerate}
The proof is complete.
\end{proof}

\subsection{Computing $\mathcal{X}^{\rm irr}_{{\rm SL}(3,\mathbb{C})}(\Gamma)$ for general $\Gamma$}

For any two-generator group $\Gamma=\langle z_1,z_2\mid r_1,\ldots,r_k\rangle$, we can decompose
\begin{align*}
\mathcal{R}^{\rm irr}_{{\rm SL}(3,\mathbb{C})}(\Gamma)&=
\mathcal{R}^{\rm asym}_{{\rm SL}(3,\mathbb{C})}(\Gamma)
\cup\mathcal{R}^{\rm irr,sym}_{{\rm SL}(3,\mathbb{C})}(\Gamma),  \\
\mathcal{X}^{\rm irr}_{{\rm SL}(3,\mathbb{C})}(\Gamma)&=
\mathcal{X}^{\rm asym}_{{\rm SL}(3,\mathbb{C})}(\Gamma)
\cup\mathcal{X}^{\rm irr,sym}_{{\rm SL}(3,\mathbb{C})}(\Gamma);
\end{align*}
the meaning of each term is self-evident.

For each $j$, the condition $\sigma_{\mathbf{a},\mathbf{b}}(r_j)=\mathbf{e}$ can be transformed into an equation of the form $\mathbf{s}_j=\mathbf{0}$, where $\mathbf{s}_j$ is a noncommutative polynomial in $\mathbf{a},\overline{\mathbf{a}}, \mathbf{b},\overline{\mathbf{b}}$.

The asymmetric part $\mathcal{X}^{\rm asym}_{{\rm SL}(3,\mathbb{C})}(\Gamma)$ is relatively easy to determine: for $j=1,\ldots,k$, $\mathbf{s}_j=\mathbf{0}$ is equivalent to $((\mathbf{u},\mathbf{s}_j))=0$ for each $\mathbf{u}\in\mathfrak{B}_\theta$ or $((\mathbf{u},\mathbf{s}_j))=0$ for each $\mathbf{u}\in\mathfrak{B}_\delta$; a system of polynomial equations in traces can be written down.

To find the symmetric part $\mathcal{X}^{\rm irr,sym}_{{\rm SL}(3,\mathbb{C})}(\Gamma)$, we may assume $\mathbf{a},\mathbf{b}$ to be symmetric. 
When $\mathbf{a}$ is special, by Corollary \ref{cor:one-is-special}, each $\mathbf{s}_j$ can be written as a linear combination of $\mathbf{e},\mathbf{a},\mathbf{b},\overline{\mathbf{b}},\mathbf{a}\mathbf{b},\mathbf{b}\mathbf{a},\mathbf{a}\overline{\mathbf{b}},\overline{\mathbf{b}}\mathbf{a},\mathbf{b}\lrcorner\mathbf{a}$ which form a basis of $\mathcal{M}$, so $\mathbf{s}_j=\mathbf{0}$ is equivalent to the simultaneous vanishing of the coefficients. When $\mathbf{b}$ is special, the situation is similar.
When $\mathbf{a}$, $\mathbf{b}$ are both ordinary,  
requiring ${\rm tr}(\mathbf{u}\mathbf{s}_j)=0$ for each 
$$\mathbf{u}\in\{\mathbf{e},\mathbf{a},\overline{\mathbf{a}},\mathbf{b},\overline{\mathbf{b}},
\mathbf{a}\mathbf{b},\mathbf{b}\mathbf{a},\overline{\mathbf{a}}\mathbf{b},\mathbf{b}\overline{\mathbf{a}},\mathbf{a}\overline{\mathbf{b}},
\overline{\mathbf{a}}\overline{\mathbf{b}}\},$$
the result is a system of polynomial equations in traces. 

When the relation $r_j$ admits a certain symmetry (this occurs for rational link groups; see \cite{BZ03} Exercise 12.1), the computation may be dramatically simplified by investigating ${\rm tr}(\mathbf{u}\mathbf{s}_j)=0$ for $\mathbf{u}$ in the bases given by Theorem \ref{thm:basis-1} and \ref{thm:basis-2}.

Explicit representations can be found without much difficulty, once the corresponding characters are determined. 

\begin{exmp}
\rm Let $\Gamma=\langle x,y\mid xy\overline{x}yx=yx\overline{y}xy\rangle$, the group of the figure eight knot.
Given $\mathbf{x},\mathbf{y}\in{\rm SL}(3,\mathbb{C})$, there exists a representation $\rho:\Gamma\to{\rm SL}(3,\mathbb{C})$ with $\rho(x)=\mathbf{x}, \rho(y)=\mathbf{y}$ if and only if
$$\mathbf{c}:=\mathbf{x}\mathbf{y}\overline{\mathbf{x}}\mathbf{y}\mathbf{x}
-\mathbf{y}\mathbf{x}\overline{\mathbf{y}}\mathbf{x}\mathbf{y}=\mathbf{0}.$$
Suppose $\rho$ is irreducible. Then $\mathbf{x},\mathbf{y}$ must be both ordinary, as they are conjugate.

In the case $t_{\overline{1}2}=t_{1\overline{2}}=s$ (the reason for considering this will soon be clear),
\begin{align*}
\mathbf{x}\mathbf{y}\overline{\mathbf{x}}\mathbf{y}\mathbf{x}
\stackrel{(\ref{eq:identity3})}=\ &\mathbf{x}(t_{\overline{1}2}\mathbf{y}-t_{1\overline{2}}\mathbf{x}+\mathbf{x}\overline{\mathbf{y}}\mathbf{x}) =t_{\overline{1}2}\mathbf{x}\mathbf{y}\mathbf{x}-t_{1\overline{2}}\mathbf{x}^3+\mathbf{x}^2\overline{\mathbf{y}}\mathbf{x}^2 \\
=\ &-\theta(\mathbf{x}\overline{\mathbf{y}})-t_{\overline{1}2}\theta(\overline{\mathbf{x}}\mathbf{y})+\overline{t}t_{\overline{1}2}\mathbf{y}
-\overline{t}^2\overline{\mathbf{y}}
+(\overline{t}-t\overline{t}^2+\overline{t}t_{1\overline{2}}+t_{12}t_{\overline{1}2}+tt_{\overline{12}})\mathbf{x} \\
&+(t_{\overline{12}}+tt_{\overline{1}2}-2\overline{t}^2)\overline{\mathbf{x}}+(1-t\overline{t})t_{\overline{1}2}\mathbf{e}+(\star),
\end{align*}
where $(\star)$ stands for the sum of terms that are symmetric between $x$ and $y$, and we have omit some intermediate steps that can be easily filled in. So
\begin{align*}
\mathbf{c}=(1-s)(\theta(\overline{\mathbf{x}}\mathbf{y})-\theta(\mathbf{x}\overline{\mathbf{y}}))
+(\overline{t}-t\overline{t}^2+st_{12}+tt_{\overline{12}})(\mathbf{x}-\mathbf{y})
+(t_{\overline{12}}+st-\overline{t}^2)(\overline{\mathbf{x}}-\overline{\mathbf{y}}).  
\end{align*}

\begin{enumerate}
  \item[(i)] If $t_{12\overline{12}}\ne t_{21\overline{21}}$, then it follows from
        \begin{align*}
        0&={\rm tr}(\mathbf{c}\delta(\overline{\mathbf{x}}\mathbf{y}))=(t_{12\overline{12}}-t_{21\overline{21}})(1-t_{1\overline{2}}), \\
        0&={\rm tr}(\mathbf{c}\delta(\overline{\mathbf{x}\mathbf{y}}))=(t_{12\overline{12}}-t_{21\overline{21}})(t_{\overline{1}2}-1)
        \end{align*}
        that $t_{\overline{1}2}=t_{1\overline{2}}=1$. Hence by (\ref{eq:cube}), $(\overline{\mathbf{x}}\mathbf{y})^4=\mathbf{e}$, and
        $$\mathbf{0}=\mathbf{c}=(\overline{t}-t\overline{t}^2+t_{12}+tt_{\overline{12}})(\mathbf{x}-\mathbf{y})
        +(t_{\overline{12}}+t-\overline{t}^2)(\overline{\mathbf{x}}-\overline{\mathbf{y}}).$$
        Similarly as in the first step of proving Theorem \ref{thm:basis-2} (ii),
        $$\overline{t}-t\overline{t}^2+t_{12}+tt_{\overline{12}}=t_{\overline{12}}+t-\overline{t}^2=0,$$
        i.e. $t_{\overline{12}}=\overline{t}^2-t$, $t_{12}=t^2-\overline{t}$.
        This implies $t_{1^2\overline{2}}=t_{\overline{1}2^2}=t_{\overline{1}^22}=t_{1\overline{2}^2}=0$, which is equivalent to $(\mathbf{x}^2\overline{\mathbf{y}})^3=(\mathbf{y}^2\overline{\mathbf{x}})^3=\mathbf{e}$.
        Thus, $\rho$ factors through
        $$\Gamma\to D(3,3,4):=\langle a,b\mid a^3=b^3=(ab)^4=1\rangle, \qquad x\mapsto bab, \ \
        y\mapsto \overline{a}\overline{b}\overline{a},$$
        which is surjective, as witnessed by $\overline{y}x\overline{y}\mapsto a$ and $x\overline{y}x\mapsto b$.
  \item[\rm(ii)] If $t_{\overline{1}2}\ne t_{1\overline{2}}$, then with $x=\overline{y}zy\overline{z}y$, the presentation of
        $\Gamma$ becomes $\langle z,y\mid zy\overline{z}yz=yz\overline{y}zy\rangle$. We have $z=x\lrcorner y$.
        Hence ${\rm tr}([\mathbf{y},\mathbf{z}])=t_{\overline{1}2}$ and ${\rm tr}([\mathbf{z},\mathbf{y}])=t_{1\overline{2}}$.
        Now that the relation has the same form as the previous one and ${\rm tr}([\mathbf{y},\mathbf{z}])\ne{\rm tr}([\mathbf{z},\mathbf{y}])$, an argument same as above leads to
        $$(\mathbf{z}^2\overline{\mathbf{y}})^3=(\mathbf{y}^2\overline{\mathbf{z}})^3=(\overline{\mathbf{z}}\mathbf{y})^4=\mathbf{e}.$$
        This time $\rho$ factors through the epimorphism
        $\Gamma\twoheadrightarrow D(3,3,4)$, $x\mapsto \overline{b}a$, $y\mapsto \overline{a}\overline{b}\overline{a}$.
  \item[\rm(iii)] Now suppose $t_{12\overline{12}}=t_{21\overline{21}}$ and $t_{\overline{1}2}=t_{1\overline{2}}=s$.
        Then $\mathbf{c}=\mathbf{0}$ is equivalent to ${\rm tr}(\mathbf{c}\mathbf{u})=0$ for $\mathbf{u}\in\{\mathbf{e},\mathbf{x},\overline{\mathbf{x}},\mathbf{y},\overline{\mathbf{y}},\theta(\mathbf{x}\mathbf{y}),
        \theta(\overline{\mathbf{x}}\mathbf{y})\}$.
        It is easy to verify
        $${\rm tr}(\mathbf{c})={\rm tr}(\mathbf{c}\theta(\mathbf{x}\mathbf{y}))=0, \qquad
        {\rm tr}(\mathbf{c}\mathbf{x})={\rm tr}(\mathbf{c}\mathbf{y}), \qquad
        {\rm tr}(\mathbf{c}\overline{\mathbf{x}})={\rm tr}(\mathbf{c}\overline{\mathbf{y}}).$$
        Thus, $\mathcal{X}^{\rm irr,sym}_{{\rm SL}(3,\mathbb{C})}(\Gamma)$ is given by
        $${\rm tr}(\mathbf{c}\mathbf{x})={\rm tr}(\mathbf{c}\overline{\mathbf{x}})={\rm tr}(\mathbf{c}\theta(\overline{\mathbf{x}}\mathbf{y}))=0,$$
        which can be written as a system of polynomial equations in $t,\overline{t},t_{12},t_{\overline{12}},s$.
\end{enumerate}

These can be compared with the components $V_0,V_1,V_2$ of $\mathcal{X}^{\rm irr}_{{\rm SL}(3,\mathbb{C})}(\Gamma)$, given in \cite{HMP16}; $S,T$ there respectively correspond to $x,y$ here. It is easy to see that $V_1,V_2$ of \cite{HMP16} are reproduced exactly by (i),(ii). We believe that $V_0$ can be identified with the result of (iii), but some more efforts are needed.
\end{exmp}

\section{Application to the double twist link groups}

\begin{figure}[h]
  \centering
  \includegraphics[width=10cm]{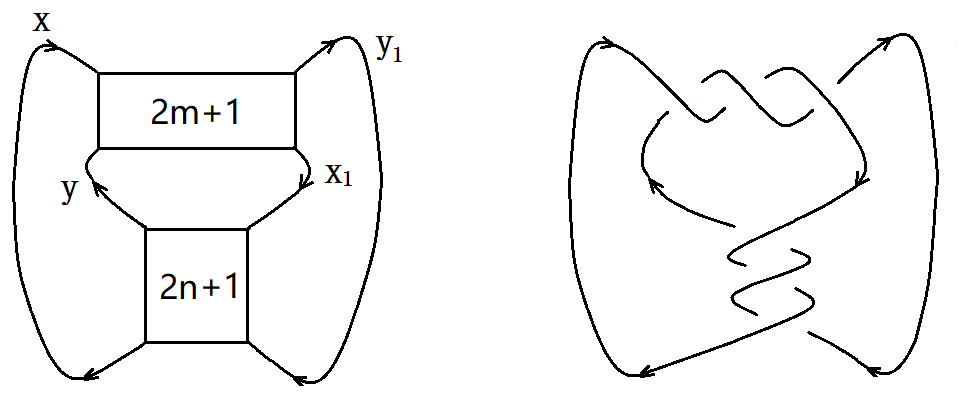}\\
  \caption{Left: the double twist link $J(2m+1,2n+1)$. Right: when $m=n=1$, it is the Whitehead link.}\label{fig:double-twist-link}
\end{figure}

Let $L=J(2m+1,2n+1)$ be the double twist link as shown in Figure \ref{fig:double-twist-link}.
Let $\Gamma=\pi_1(E_L)$, where $E_L$ is the complement of a tubular neighborhood of $L$ in $S^3$.
Adopting Wirtinger presentation, we have $x_1=z^{m}\lrcorner x$ and $y_1=z^{m+1}\lrcorner y$, where $z=xy$;
writing $y_1$ in another way leads to the unique relation
$(x_1\overline{y})^{n+1}\lrcorner y=z^{m+1}\lrcorner y.$
Substituting $x=z\overline{y}$ results in
\begin{align*}
\Gamma=\langle y,z\mid \overline{y}z^{m+1}w^{n+1}=z^{m+1}w^{n+1}\overline{y}\rangle, \qquad \text{with} \quad w=\overline{z}^{m+1}yz^{m}y.
\end{align*}
Note that $z^{m+1}w^{n+1}=(z^{m+1}w^{n+1})^\leftarrow$, so the relation is the same as $$\overline{y}z^{m+1}w^{n+1}=(\overline{y}z^{m+1}w^{n+1})^\leftarrow.$$

Given $\mathbf{y},\mathbf{z}\in{\rm SL}(3,\mathbb{C})$, there exists a representation $\rho:\Gamma\to{\rm SL}(3,\mathbb{C})$ with $\rho(y)=\mathbf{y},\rho(z)=\mathbf{z}$ if and only if
\begin{align}
\mathbf{c}:=\delta(\overline{\mathbf{y}}\mathbf{z}^{m+1}\mathbf{w}^{n+1})=\mathbf{0}, \qquad \text{with} \quad
\mathbf{w}=\overline{\mathbf{z}}^{m+1}\mathbf{y}\mathbf{z}^{m}\mathbf{y}. \label{eq:c-origin}
\end{align}

Abbreviate $\alpha_k(\mathbf{z}),\beta_k(\mathbf{z}),\gamma_k(\mathbf{z})$ to $\alpha_k,\beta_k,\gamma_k$, respectively.
Let $\tilde{\alpha}=\alpha_{n+1}(\mathbf{w})$, $\tilde{\beta}=\beta_{n+1}(\mathbf{w})$, $\tilde{\gamma}=\gamma_{n+1}(\mathbf{w})$.
Then
\begin{align}
\mathbf{c}&=\delta(\overline{\mathbf{y}}\mathbf{z}^{m+1}(\tilde{\alpha}\mathbf{w}+\tilde{\beta}\mathbf{e}+\tilde{\gamma}\overline{\mathbf{w}})) \nonumber  \\
&=\tilde{\alpha}\delta(\mathbf{z}^{m}\mathbf{y})+\tilde{\beta}\delta(\overline{\mathbf{y}}\mathbf{z}^{m+1})+
\tilde{\gamma}\delta(\overline{\mathbf{y}}\mathbf{z}^{m+1}\overline{\mathbf{y}}\overline{\mathbf{z}}^m\overline{\mathbf{y}}\mathbf{z}^{m+1}).    \label{eq:c-1}
\end{align}

To shorten the writing, let $p_k={\rm tr}(\mathbf{y}\mathbf{z}^k)$, $q_k={\rm tr}(\overline{\mathbf{y}}\mathbf{z}^k)$.
Then
\begin{align*}
\overline{\mathbf{y}}\mathbf{z}^{m+1}\overline{\mathbf{y}}\overline{\mathbf{z}}^m\overline{\mathbf{y}}\mathbf{z}^{m+1}&=\overline{\mathbf{y}}\mathbf{z}^{m+1}\big(q_{-m}\overline{\mathbf{y}}-p_{m}\mathbf{z}^{m}
+\mathbf{z}^{m}\mathbf{y}\mathbf{z}^{m}\big)\mathbf{z}^{m+1} \\
&=q_{-m}(\overline{\mathbf{y}}\mathbf{z}^{m+1})^2-p_{m}\overline{\mathbf{y}}\mathbf{z}^{3m+2}+\overline{\mathbf{y}}\mathbf{z}^{2m+1}\mathbf{y}\mathbf{z}^{2m+1} \\
&=q_{-m}\big(q_{m+1}\overline{\mathbf{y}}\mathbf{z}^{m+1}-p_{-m-1}\mathbf{e}+\overline{\mathbf{z}}^{m+1}\mathbf{y}\big)-p_{m}\overline{\mathbf{y}}\mathbf{z}^{3m+2} \\
&\ \ \ +\overline{\mathbf{y}}(p_{2m+1}\mathbf{z}^{2m+1}-q_{-2m-1}\overline{\mathbf{y}}+\overline{\mathbf{y}}\overline{\mathbf{z}}^{2m+1}\overline{\mathbf{y}}).
\end{align*}
Hence
\begin{align*}
&\delta(\overline{\mathbf{y}}\mathbf{z}^{m+1}\overline{\mathbf{y}}\overline{\mathbf{z}}^m\overline{\mathbf{y}}\mathbf{z}^{m+1}) \\
=\ &q_{-m}q_{m+1}\delta(\overline{\mathbf{y}}\mathbf{z}^{m+1})+q_{-m}\delta(\overline{\mathbf{z}}^{m+1}\mathbf{y})
-p_{m}\delta(\overline{\mathbf{y}}\mathbf{z}^{3m+2})+p_{2m+1}\delta(\overline{\mathbf{y}}\mathbf{z}^{2m+1}) \\
&+p_0\delta(\overline{\mathbf{y}}\overline{\mathbf{z}}^{2m+1})+\delta(\mathbf{y}\overline{\mathbf{z}}^{2m+1}\overline{\mathbf{y}}) \\
=\ &q_{-m}q_{m+1}(\alpha_{m+1}\delta(\overline{\mathbf{y}}\mathbf{z})+\gamma_{m+1}\delta(\overline{\mathbf{y}}\overline{\mathbf{z}}))
+q_{-m}(\alpha_{-m-1}\delta(\mathbf{z}\mathbf{y})+\gamma_{-m-1}\delta(\overline{\mathbf{z}}\mathbf{y})) \\
&-p_{m}(\alpha_{3m+2}\delta(\overline{\mathbf{y}}\mathbf{z})+\gamma_{3m+2}\delta(\overline{\mathbf{y}}\overline{\mathbf{z}}))
+p_{2m+1}(\alpha_{2m+1}\delta(\overline{\mathbf{y}}\mathbf{z})+\gamma_{2m+1}\delta(\overline{\mathbf{y}}\overline{\mathbf{z}})) \\
&+p_0(\alpha_{-2m-1}\delta(\overline{\mathbf{y}}\mathbf{z})+\gamma_{-2m-1}\delta(\overline{\mathbf{y}}\overline{\mathbf{z}}))
+\alpha_{-2m-1}\delta(\mathbf{y}\lrcorner\mathbf{z})+\gamma_{-2m-1}\delta(\mathbf{y}\lrcorner\overline{\mathbf{z}}) \\
=&-q_{-m}\alpha_{-m-1}\delta(\mathbf{y}\mathbf{z})+\varsigma_1\delta(\overline{\mathbf{y}}\mathbf{z})
-q_{-m}\gamma_{-m-1}\delta(\mathbf{y}\overline{\mathbf{z}})+\varsigma_2\delta(\overline{\mathbf{y}}\overline{\mathbf{z}}) \\
&+\alpha_{-2m-1}\delta(\mathbf{y}\lrcorner\mathbf{z})+\gamma_{-2m-1}\delta(\mathbf{y}\lrcorner\overline{\mathbf{z}}),
\end{align*}
with
\begin{align}
\varsigma_1&=q_{-m}q_{m+1}\alpha_{m+1}-p_m\alpha_{3m+2}+p_{2m+1}\alpha_{2m+1}+p_0\alpha_{-2m-1},  \label{eq:varsigma1} \\
\varsigma_2&=q_{-m}q_{m+1}\gamma_{m+1}-p_m\gamma_{3m+2}+p_{2m+1}\gamma_{2m+1}+p_0\gamma_{-2m-1}.  \label{eq:varsigma2}
\end{align}

By (\ref{eq:c-1}),
\begin{align}
\mathbf{c}=\eta_1\delta(\mathbf{y}\mathbf{z})+\eta_2\delta(\overline{\mathbf{y}}\mathbf{z})+\eta_3\delta(\mathbf{y}\overline{\mathbf{z}})
+\eta_4\delta(\overline{\mathbf{y}\mathbf{z}})+\eta_5\delta(\mathbf{y}\lrcorner\mathbf{z})+\eta_6\delta(\mathbf{y}\lrcorner\overline{\mathbf{z}}), \label{eq:c-2}
\end{align}
with
\begin{align}
\eta_1&=-\tilde{\gamma}q_{-m}\alpha_{-m-1}-\tilde{\alpha}\alpha_{m}, && \eta_2=\tilde{\gamma}\varsigma_1+\tilde{\beta}\alpha_{m+1}, \label{eq:eta1-2} \\
\eta_3&=-\tilde{\gamma}q_{-m}\gamma_{-m-1}-\tilde{\alpha}\gamma_{m}, &&\eta_4=\tilde{\gamma}\varsigma_2+\tilde{\beta}\gamma_{m+1}, \label{eq:eta3-4} \\
\eta_5&=\tilde{\gamma}\alpha_{-2m-1}, &&\eta_6=\tilde{\gamma}\gamma_{-2m-1}.  \label{eq:eta5-6}
\end{align}

\subsection{The asymmetric slice: ${\rm tr}(\mathbf{y}\mathbf{z}\overline{\mathbf{y}\mathbf{z}})\ne{\rm tr}(\mathbf{z}\mathbf{y}\overline{\mathbf{z}\mathbf{y}})$}

The table given in Section \ref{sec:assymetric} will be repeatedly applied.

By Theorem \ref{thm:basis-1}, $\mathbf{c}=\mathbf{0}$ if and only if
$$((\mathbf{u},\mathbf{c}))=0, \qquad \text{for} \quad  \mathbf{u}\in\{\mathbf{e},\mathbf{y},\overline{\mathbf{y}},\mathbf{z},\overline{\mathbf{z}},\theta(\mathbf{y}\mathbf{z}),\theta(\mathbf{y}\overline{\mathbf{z}}),
\theta(\overline{\mathbf{y}}\mathbf{z}),\theta(\overline{\mathbf{y}\mathbf{z}})\}.$$
From (\ref{eq:c-2}) it is clear that $((\mathbf{u},\mathbf{c}))=0$ for $\mathbf{u}\in\{\mathbf{e},\mathbf{y},\overline{\mathbf{y}}\}$.
As necessary conditions,
$0=((\mathbf{z},\mathbf{c}))=-\tilde{\gamma}\gamma_{-2m-1}$, and $0=((\overline{\mathbf{z}},\mathbf{c}))=\tilde{\gamma}\alpha_{-2m-1}.$

\begin{itemize}
  \item If $\tilde{\gamma}\ne 0$, then $\alpha_{-2m-1}=\gamma_{-2m-1}=0$, implying $\mathbf{z}^{2m+1}=\mu\mathbf{e}$ with $\mu\in\langle\omega\rangle$, so $\mathbf{w}=\mu^{-1}(\mathbf{z}^m\mathbf{y})^{2}$, and (\ref{eq:c-origin}) becomes
      $$\mathbf{0}=\mathbf{c}':=\mu^{n}\mathbf{c}=\delta((\mathbf{z}^m\mathbf{y})^{2n+1}).$$
      Let $\alpha'=\alpha_{2n+1}(\mathbf{z}^m\mathbf{y})$, $\gamma'=\gamma_{2n+1}(\mathbf{z}^m\mathbf{y})$.
      Then
      \begin{align*}
      \mathbf{c}'&=\alpha'\delta(\mathbf{z}^m\mathbf{y})+\gamma'\delta(\overline{\mathbf{y}}\overline{\mathbf{z}}^m) \\
      &=\alpha'(\alpha_m\delta(\mathbf{z}\mathbf{y})+\gamma_m\delta(\overline{\mathbf{z}}\mathbf{y}))
      +\gamma'(\alpha_{-m}\delta(\overline{\mathbf{y}}\mathbf{z})+\gamma_{-m}\delta(\overline{\mathbf{y}}\overline{\mathbf{z}})).
      \end{align*}
      We easily obtain
      \begin{align*}
      0&=((\theta(\mathbf{y}\mathbf{z}),\mathbf{c}'))=\gamma_{-m}\gamma',
      && 0=((\theta(\mathbf{y}\overline{\mathbf{z}}),\mathbf{c}'))=-\alpha_{-m}\gamma', \\
      0&=((\theta(\overline{\mathbf{y}}\mathbf{z}),\mathbf{c}'))=\gamma_m\alpha',
      &&0=((\theta(\overline{\mathbf{y}\mathbf{z}}),\mathbf{c}'))=-\alpha_{m}\alpha'.
      \end{align*}
      If $\alpha_m=\gamma_m=0$, then $\mathbf{z}^m=\kappa\mathbf{e}$, $\kappa\in\langle\omega\rangle$, but this together with $\mathbf{z}^{2m+1}=\mu\mathbf{e}$ would imply $\mathbf{z}=\kappa\mu\mathbf{e}$. Similarly, $\alpha_{-m}=\gamma_{-m}=0$ is neither possible.
      Hence $\alpha'=\gamma'=0$, implying $(\mathbf{z}^m\mathbf{y})^{2n+1}=\nu\mathbf{e}$, $\nu\in\langle\omega\rangle$.
  \item If $\tilde{\gamma}=0$, then
        \begin{align*}
        \mathbf{c}&=\tilde{\alpha}\delta(\mathbf{z}^{m}\mathbf{y})+\tilde{\beta}\delta(\overline{\mathbf{y}}\mathbf{z}^{m+1}) \\
        &=\tilde{\alpha}(\alpha_m\delta(\mathbf{z}\mathbf{y})+\gamma_m\delta(\overline{\mathbf{z}}\mathbf{y}))
        +\tilde{\beta}(\alpha_{m+1}\delta(\overline{\mathbf{y}}\mathbf{z})+\gamma_{m+1}\delta(\overline{\mathbf{y}}\overline{\mathbf{z}})).
        \end{align*}
        As is easy to verify, $((\mathbf{z},\mathbf{c}))=((\overline{\mathbf{z}},\mathbf{c}))=0$, and
        \begin{align*}
        0&=((\theta(\mathbf{y}\mathbf{z}),\mathbf{c}))=\gamma_{m+1}\tilde{\beta},
        &&0=((\theta(\mathbf{y}\overline{\mathbf{z}}),\mathbf{c}))=-\alpha_{m+1}\tilde{\beta}, \\
        0&=((\theta(\overline{\mathbf{y}}\mathbf{z}),\mathbf{c}))=\gamma_{m}\tilde{\alpha},
        &&0=((\theta(\overline{\mathbf{y}\mathbf{z}}),\mathbf{c}))=-\alpha_{m}\tilde{\alpha}.
        \end{align*}
        There are two possibilities.
        \begin{itemize}
          \item If $\alpha_{m}=\gamma_{m}=0$, then $\mathbf{z}^m=\mu\mathbf{e}$, $\mu\in\langle\omega\rangle$. Consequently, $\alpha_{m+1}\ne 0$ or $\gamma_{m+1}\ne 0$, either implying $\tilde{\beta}=0$. Then $\mathbf{w}^{n+1}=\tilde{\alpha}\mathbf{w}$, which can be rewritten as $\mathbf{w}^n=\nu\mathbf{e}$, $\nu\in\langle\omega\rangle$; equivalently, $(\overline{\mathbf{z}}\mathbf{y}^2)^{n}=\nu\mathbf{e}$.
          \item Otherwise, $\tilde{\alpha}=0$, so that $\mathbf{w}^{n+1}=\nu\mathbf{e}$, $\nu\in\langle\omega\rangle$, and $\tilde{\beta}=\nu$. Hence
              $\alpha_{m+1}=\gamma_{m+1}=0$, implying $\mathbf{z}^{m+1}=\mu\mathbf{e}$, $\mu\in\langle\omega\rangle$. Then $\mathbf{w}^{n+1}=\nu\mathbf{e}$ becomes $(\overline{\mathbf{z}}\mathbf{y}^2)^{n+1}=\nu\mathbf{e}$.
        \end{itemize}
\end{itemize}

\begin{lem}\label{lem:power}
Let $k$ be a nonzero integer. If $\mathbf{z}$ is ordinary and $\mathbf{z}^k=\mu\mathbf{e}$ with $\mu\in\langle\omega\rangle$, then
$|k|\ge3$, $\mathbf{z}$ is conjugate to $\mathbf{d}^{\lambda}_{\lambda\zeta}$ for some $\lambda,\zeta$ such that
$\lambda^k=\mu$, $\zeta^k=1$, and $\zeta\notin\{1,\lambda^{-3},\pm\sqrt{\lambda^{-3}}\}.$
In particular, $\mu=1$ when $k=\pm3$.
\end{lem}

\begin{proof}
Now that $\mathbf{z}^k=\mu\mathbf{e}$ and $\mathbf{z}$ is ordinary, $|k|$ must be at least 3.
It is impossible for $\mathbf{z}$ to be conjugate to $\mathbf{h}_\lambda$ or $\nu\mathbf{k}$, as
\begin{align*}
(\mathbf{h}_\lambda)^k&=(\mathbf{d}^\lambda_\lambda+\lambda\mathbf{f}^2)^k=(\mathbf{d}^\lambda_\lambda)^k+k\lambda(\mathbf{d}^\lambda_\lambda)^{k-1}\mathbf{f}^2\notin\mathbb{C}\langle\mathbf{e}\rangle, \\
\mathbf{k}^k&=(\mathbf{e}+\mathbf{f}+\mathbf{f}^2)^k=\mathbf{e}+k\mathbf{f}+\frac{k(k+1)}{2}\mathbf{f}^2\notin\mathbb{C}\langle\mathbf{e}\rangle.
\end{align*}
Hence $\mathbf{z}$ is conjugate to $\mathbf{d}^{\lambda}_{\lambda\zeta}$ for some $\lambda,\zeta$ such that $\lambda^k=\mu$, $\zeta^k=1$. Since $\mathbf{z}$ is ordinary, $\lambda,\lambda\zeta,(\lambda^2\zeta)^{-1}$ should be distinct; equivalently, $\zeta\notin\{1,\lambda^{-3},\pm\sqrt{\lambda^{-3}}\}$. This constraint forces $\mu=1$ when $k=\pm3$.
\end{proof}

Put $\Lambda_k=\emptyset$ if $|k|\le 2$, and put
$$\Lambda_k=\{(\lambda,\zeta)\colon\lambda^{3k}=\zeta^k=1,\ \zeta\ne 1,\lambda^{-3},\pm\sqrt{\lambda^{-3}}\}$$
if $|k|\ge 3$. Recall (\ref{eq:char-variety-F2}):
$$\mathcal{X}_{{\rm SL}(3,\mathbb{C})}(F_2)\cong\{(s_1,s_{\overline{1}},s_2,s_{\overline{2}},s_3,s_{\overline{3}}, s_4, s_{\overline{4}},s_5)\colon s_5^2-Ps_5+Q=0\}.$$

\begin{thm} \label{thm:asym}
The asymmetric slice $\mathcal{X}^{\rm asym}_{{\rm SL}(3,\mathbb{C})}(\Gamma)$ is the disjoint union of
\begin{align*}
&\mathcal{X}^1_{\lambda_1,\zeta_1;\lambda_2,\zeta_2}, \qquad (\lambda_1,\zeta_1)\in\Lambda_{2m+1},  \quad (\lambda_2,\zeta_2)\in\Lambda_{2n+1},  \\
&\mathcal{X}^2_{\lambda_1,\zeta_1;\lambda_2,\zeta_2}, \qquad  (\lambda_1,\zeta_1)\in\Lambda_{m},  \qquad  (\lambda_2,\zeta_2)\in\Lambda_{n},  \\
&\mathcal{X}^3_{\lambda_1,\zeta_1;\lambda_2,\zeta_2}, \qquad (\lambda_1,\zeta_1)\in\Lambda_{m+1}, \quad (\lambda_2,\zeta_2)\in\Lambda_{n+1},
\end{align*}
where $\mathcal{X}^j_{\lambda_1,\zeta_1;\lambda_2,\zeta_2}$ is isomorphic to the subspace of $\mathcal{X}_{{\rm SL}(3,\mathbb{C})}(F_2)$ determined by
$$s_i=\lambda_i+\lambda_i\zeta_i+(\lambda_i^2\zeta_i)^{-1}, \ \ \ s_{\overline{i}}=\lambda_i^{-1}+(\lambda_i\zeta_i)^{-1}+\lambda_i^2\zeta_i, \ \ \ i=1,2, \qquad P^2\ne 4Q.$$
\end{thm}

\begin{proof}
According to the above discussions, there are three possibilities for $\mathbf{y},\mathbf{z}$ to define a representation:
1) $\mathbf{z}^{2m+1}=\mu\mathbf{e}$, $(\mathbf{z}^m\mathbf{y})^{2n+1}=\nu\mathbf{e}$;
2) $\mathbf{z}^m=\mu\mathbf{e}$, $(\overline{\mathbf{z}}\mathbf{y}^2)^{n}=\nu\mathbf{e}$;
3) $\mathbf{z}^{m+1}=\mu\mathbf{e}$, $(\overline{\mathbf{z}}\mathbf{y}^2)^{n+1}=\nu\mathbf{e}$.
In each case, $\mu,\nu\in\langle\omega\rangle$.

In case 1), let $\mathbf{a}=\mathbf{z}$, $\mathbf{b}=\mathbf{z}^m\mathbf{y}$. 
Then ${\rm tr}(\mathbf{y}\mathbf{z}\overline{\mathbf{y}\mathbf{z}})\ne{\rm tr}(\mathbf{z}\mathbf{y}\overline{\mathbf{z}\mathbf{y}})$ is equivalent to 
${\rm tr}(\mathbf{a}\mathbf{b}\overline{\mathbf{a}\mathbf{b}})\ne{\rm tr}(\mathbf{b}\mathbf{a}\overline{\mathbf{b}\mathbf{a}})$.
By Lemma \ref{lem:power}, $\mathbf{a}$ is conjugate to $\mathbf{d}^{\lambda_1}_{\lambda_1\zeta_1}$ for some $(\lambda_1,\zeta_1)\in\Lambda_{2m+1}$, and $\mathbf{b}$ is conjugate to $\mathbf{d}^{\lambda_2}_{\lambda_2\zeta_2}$ for some $(\lambda_2,\zeta_2)\in\Lambda_{2n+1}$. Hence 
\begin{align*}
&{\rm tr}(\mathbf{a})=\lambda_1+\lambda_1\zeta_1+(\lambda_1^2\zeta_1)^{-1}, &&{\rm tr}(\overline{\mathbf{a}})=\lambda_1^{-1}+(\lambda_1\zeta_1)^{-1}+\lambda_1^2\zeta_1, \\
&{\rm tr}(\mathbf{b})=\lambda_2+\lambda_2\zeta_2+(\lambda_2^2\zeta_2)^{-1}, &&{\rm tr}(\overline{\mathbf{b}})=\lambda_2^{-1}+(\lambda_2\zeta_2)^{-1}+\lambda_2^2\zeta_2.
\end{align*}

Cases 2) and 3) are dealt with similarly.
\end{proof}

\subsection{The symmetric slice: ${\rm tr}(\mathbf{y}\mathbf{z}\overline{\mathbf{y}\mathbf{z}})
={\rm tr}(\mathbf{z}\mathbf{y}\overline{\mathbf{z}\mathbf{y}})$}

Suppose $\mathbf{y},\mathbf{z}$ are symmetric and $(\mathbf{y},\mathbf{z})$ is irreducible. Let $t_1={\rm tr}(\mathbf{y})$, $t_2={\rm tr}(\mathbf{z})$, and so forth.
Recall (\ref{eq:varsigma1})--(\ref{eq:eta5-6}).

\begin{thm} \label{thm:sym}
The symmetric slice $\mathcal{X}^{\rm sym,irr}_{{\rm SL}(3,\mathbb{C})}(\Gamma)=\mathcal{X}_{\rm sp,od}\cup\mathcal{X}_{\rm od,sp}\cup\mathcal{X}_{\rm od,od}$, where
$\mathcal{X}_{\rm sp,od}$ is determined by
\begin{align}
\lambda\in\mathbb{C}^\ast, \qquad t_1=2\lambda+\lambda^{-2}, \qquad  t_{\overline{1}}=2\lambda^{-1}+\lambda^2, \nonumber \\  
t_{\overline{1}2}=(\lambda^2+\lambda^{-1})t_2-\lambda t_{12}, \quad t_{\overline{12}}=(\lambda^2+\lambda^{-1})t_{\overline{2}}-\lambda t_{1\overline{2}}, \nonumber  \\
\eta_1-\lambda\eta_2+(\lambda^2+\lambda^{-1})\eta_5=\eta_3-\lambda\eta_4+(\lambda^2+\lambda^{-1})\eta_6=0, \label{eq:y-special}  
\end{align}
$\mathcal{X}_{\rm od,sp}$ is determined by 
\begin{align}
\lambda \in\mathbb{C}^\ast, \qquad t_2=2\lambda+\lambda^{-2}, \qquad t_{\overline{2}}=2\lambda^{-1}+\lambda^2,  \nonumber \\
t_{1\overline{2}}=(\lambda^2+\lambda^{-1})t_1-\lambda t_{12}, \quad t_{\overline{12}}=(\lambda^2+\lambda^{-1})t_{\overline{1}}-\lambda t_{\overline{1}2}, \nonumber  \\
\eta_1-\lambda\eta_3=\eta_2-\lambda\eta_4=\eta_5-\lambda\eta_6=0, \label{eq:z-special}
\end{align}
and $\mathcal{X}_{\rm od,od}$ is determined by $[\eta_1,\ldots,\eta_6]A=0$, with
\begin{align}
A=\left[\begin{array}{cccccc}
t_{1212}-t_{1^22^2} & t_{12\overline{1}2}-t_{2^2} & t_{121\overline{2}}-t_{1^2} &  t_{12\overline{12}}-3  \\
t_{12\overline{1}2}-t_{2^2}  & t_{\overline{1}2\overline{1}2}-t_{\overline{1}^22^2} & t_{12\overline{12}}-3 & t_{\overline{1}2\overline{12}}-t_{\overline{1}^2} \\
t_{121\overline{2}}-t_{1^2} & t_{12\overline{12}}-3 & t_{1\overline{2}1\overline{2}}-t_{1^2\overline{2}^2} &  t_{1\overline{212}}-t_{\overline{2}^2} \\
t_{12\overline{12}}-3 & t_{\overline{1}2\overline{12}}-t_{\overline{1}^2} &t_{1\overline{212}}-t_{\overline{2}^2}  & t_{\overline{1212}}-t_{\overline{1}^2\overline{2}^2}  \\
t_{12^2}-t_{1^22\overline{1}2}  & t_{12\overline{1}^22}-t_{\overline{1}2^2} & t_1 -t_{1^22\overline{12}} &  t_{12\overline{1}^2\overline{2}}-t_{\overline{1}} \\
t_1-t_{1^2\overline{21}2} &  t_{1\overline{2}\overline{1}^22}-t_{\overline{1}} & t_{1\overline{2}^2}-t_{1^2\overline{212}} & t_{1\overline{2}\overline{1}^2\overline{2}}-t_{\overline{1}\overline{2}^2}
\end{array}\right]. \label{eq:A}
\end{align}
\end{thm}

\begin{proof}
If $\mathbf{y}$ is special, with multiple eigenvalue $\lambda$, so that $\overline{\mathbf{y}}=(\lambda^2+\lambda^{-1})\mathbf{e}-\lambda\mathbf{y}$, then
$\delta(\overline{\mathbf{y}}\mathbf{z}^{\pm 1})=-\lambda\delta(\mathbf{y}\mathbf{z}^{\pm 1})$, and $\delta(\mathbf{y}\lrcorner\mathbf{z}^{\pm1})=(\lambda^2+\lambda^{-1})\delta(\mathbf{y}\mathbf{z}^{\pm1}),$
implying
$$\mathbf{c}=(\eta_1-\lambda\eta_2+(\lambda^2+\lambda^{-1})\eta_5)\delta(\mathbf{y}\mathbf{z})+(\eta_3-\lambda\eta_4+(\lambda^2+\lambda^{-1})\eta_6)\delta(\mathbf{y}\overline{\mathbf{z}}).$$
By Theorem \ref{thm:basis-2} (i), $\mathbf{c}=\mathbf{0}$ is equivalent to (\ref{eq:y-special}).

If $\mathbf{z}$ is special, with multiple eigenvalue $\lambda$, so that $\overline{\mathbf{z}}=(\lambda^2+\lambda^{-1})\mathbf{e}-\lambda\mathbf{z}$, then
$\delta(\mathbf{y}^{\pm 1}\overline{\mathbf{z}})=-\lambda\delta(\mathbf{y}^{\pm 1}\mathbf{z})$, and $\delta(\mathbf{y}\lrcorner\overline{\mathbf{z}})=-\lambda\delta(\mathbf{y}\lrcorner\mathbf{z}),$
implying
$$\mathbf{c}=(\eta_1-\lambda\eta_3)\delta(\mathbf{y}\mathbf{z})+(\lambda^2-\lambda\eta_4)\delta(\overline{\mathbf{y}}\mathbf{z})+(\eta_5-\lambda\eta_6)\delta(\mathbf{y}\lrcorner\mathbf{z}).$$
Hence $\mathbf{c}=\mathbf{0}$ is equivalent to (\ref{eq:z-special}).

Now suppose $\mathbf{y}$, $\mathbf{z}$ are both ordinary.
Let $\mathbf{u}_1=\mathbf{v}_1=\mathbf{y}\mathbf{z}$, $\mathbf{u}_2=\mathbf{v}_2=\overline{\mathbf{y}}\mathbf{z}$, $\mathbf{u}_3=\mathbf{v}_3=\mathbf{y}\overline{\mathbf{z}}$,  $\mathbf{u}_4=\mathbf{v}_4=\overline{\mathbf{y}}\overline{\mathbf{z}}$, $\mathbf{u}_5=\mathbf{y}\lrcorner\mathbf{z}$, $\mathbf{u}_6=\mathbf{y}\lrcorner\overline{\mathbf{z}}$.
Direct computation shows $[{\rm tr}(\delta(\mathbf{u}_i)\delta(\mathbf{v}_j))/2]_{6\times 4}=A$ as given by (\ref{eq:A}). 
By Theorem \ref{thm:basis-2} (ii), $\mathbf{c}=\mathbf{0}$ is equivalent to $[\eta_1,\ldots,\eta_6]A=0.$
\end{proof}

\subsection{The Whitehead link}

When $m=n=1$, $L$ is the Whitehead link. Now
$$\Gamma=\langle y,z\mid yzy\overline{z}^2yz=zy\overline{z}^2yzy\rangle.$$
With $a=z,b=zy$, one has an alternative presentation 
$\langle a,b\mid \overline{a}b^2\overline{a}^3b=b\overline{a}^3b^2\overline{a}\rangle,$
which is the same as given by (3) in \cite{GW18}. Let $\mathbf{a}=\mathbf{z},\mathbf{b}=\mathbf{z}\mathbf{y}$, and let $t_1={\rm tr}(\mathbf{a}),t_2={\rm tr}(\mathbf{b})$ and so forth.

For the asymmetric slice, by Lemma \ref{lem:power}, only case 1) in the proof of Theorem \ref{thm:asym} occurs: $\mathbf{a}^3=\mathbf{b}^3=\mathbf{e}$, which is equivalent to $t_1=t_{\overline{1}}=t_2=t_{\overline{2}}=0$. 
This recovers the result of \cite{GW18}.
With (\ref{eq:char-variety-F2}) recalled, $\mathcal{X}^{\rm asym}_{{\rm SL}(3,\mathbb{C})}(\Gamma)$ is identified with the subspace of $\mathcal{X}_{{\rm SL}(3,\mathbb{C})}(F_2)$ given by 
$$\big\{(0,0,0,0,s_3,s_{\overline{3}},s_4,s_{\overline{4}},s_5)\colon s_5^2-Ps_5+Q=0,\ P^2\ne 4Q\big\}.$$

For the symmetric slice, let $\check{\mathbf{c}}=\delta(\overline{\mathbf{a}}\mathbf{b}^2\overline{\mathbf{a}}^3\mathbf{b})$. Since
$$\overline{\mathbf{a}}\mathbf{b}^2\overline{\mathbf{a}}^3\mathbf{b}=(t_{\overline{1}}^2-t_1)\overline{\mathbf{a}}\mathbf{b}^2\overline{\mathbf{a}}\mathbf{b}+(1-t_1t_{\overline{1}})\overline{\mathbf{a}}\mathbf{b}^3+
t_{\overline{1}}\overline{\mathbf{a}}\mathbf{b}^2\mathbf{a}\mathbf{b},$$
we have
\begin{align*}
\check{\mathbf{c}}=\ &(t_{\overline{1}}^2-t_1)\big(t_{\overline{1}2^2}\delta(\overline{\mathbf{a}}\mathbf{b})+\delta(\overline{\mathbf{b}}^2\mathbf{a}\overline{\mathbf{b}})\big)
+(1-t_1t_{\overline{1}})\big((t_2^2-t_{\overline{2}})\delta(\overline{\mathbf{a}}\mathbf{b})+t_2\delta(\overline{\mathbf{a}}\overline{\mathbf{b}})\big) \\
&+t_{\overline{1}}\big(t_2\delta(\overline{\mathbf{a}}\mathbf{b}\mathbf{a}\mathbf{b})+\delta(\overline{\mathbf{a}\mathbf{b}}\mathbf{a}\mathbf{b})\big) \\
=\ &((t_{\overline{1}}^2-t_1)t_{\overline{1}2^2}+(1-t_1t_{\overline{1}})(t_2^2-t_{\overline{2}}))\delta(\overline{\mathbf{a}}\mathbf{b})+(1-t_1t_{\overline{1}})t_2\delta(\overline{\mathbf{a}}\overline{\mathbf{b}})
+t_{\overline{1}}\delta(\overline{\mathbf{a}\mathbf{b}}\mathbf{a}\mathbf{b}) \\
&+(t_{\overline{1}}^2-t_1)\big(-t_2\delta(\mathbf{a}\overline{\mathbf{b}})+\delta(\mathbf{b}\mathbf{a}\overline{\mathbf{b}})\big)+t_{\overline{1}}t_2\big(t_{12}\delta(\overline{\mathbf{a}}\mathbf{b})
-t_1\delta(\overline{\mathbf{b}}\overline{\mathbf{a}})+\delta(\mathbf{a}\overline{\mathbf{b}}\overline{\mathbf{a}})\big) \\
=\ &\eta\delta(\overline{\mathbf{a}}\mathbf{b})+t_2(t_1-t_{\overline{1}}^2)\delta(\mathbf{a}\overline{\mathbf{b}})+t_2\delta(\overline{\mathbf{a}}\overline{\mathbf{b}})
+t_{\overline{1}}t_2\delta(\mathbf{a}\lrcorner\overline{\mathbf{b}})+(t_{\overline{1}}^2-t_1)\delta(\mathbf{b}\lrcorner\mathbf{a})-t_{\overline{1}}\delta([\mathbf{b},\mathbf{a}]),
\end{align*}
with
$$\eta=(t_{\overline{1}}^2-t_1)(t_2t_{\overline{1}2}-t_{\overline{1}}t_{\overline{2}}+t_{\overline{12}})+(1-t_1t_{\overline{1}})(t_2^2-t_{\overline{2}})+t_{\overline{1}}t_2t_{12}.$$
In the last step, we have expanded $t_{\overline{1}2^2}$ as $t_2t_{\overline{1}2}-t_{\overline{1}}t_{\overline{2}}+t_{\overline{12}}$.
Similarly as in Theorem \ref{thm:sym}, we can decompose $\mathcal{X}^{\rm sym,irr}_{{\rm SL}(3,\mathbb{C})}(\Gamma)$ into three parts, corresponding respectively to three cases: $\mathbf{a}$ is special; $\mathbf{b}$ is special; $\mathbf{a}$, $\mathbf{b}$ are both ordinary.

\end{document}